\documentclass{article}
  \usepackage{amsfonts,amsmath,amssymb,amsopn,amsthm,xcolor}
  \usepackage[square,sort&compress,comma,numbers]{natbib}
  \usepackage{framed}
  \usepackage{soul}

\newcommand{\thistime}{\expandafter\calctimeA\pdfcreationdate\@nil}
\def\calctimeA#1:#2#3#4#5#6#7#8#9{\calctimeB}
\def\calctimeB#1#2#3#4#5\@nil{#1#2:#3#4}

 \usepackage[colorlinks=true]{hyperref}

\newcommand{\beq}{\begin{equation}}
\newcommand{\eeq}{\end{equation}}
  
   \def\C{\mathbb{C}}
   \def\R{\mathbb{R}}
   \def\N{\mathbb{N}}
   
   \def\1{{\rm I\mskip -10.5mu 1}}
   
   \def\e{{\varepsilon}}
   \def\D{{\nabla}}
   \def\lo{\mathop{\longrightarrow}}

   \def\cA{{\cal A}}
   
   \def\cC{{\cal C}}

   \def\cG{{\cal G}}
   \def\cH{{\cal H}}

   \def\cL{{\cal L}}
   
   \def\cN{{\cal N}}

   \def\Id{\mathop{\rm Id}\nolimits}
   \def\interior{\mathop{\rm int}\nolimits}
   
   \def\loc{\mathop{\rm loc}\nolimits}
   \def\meas{\mathop{\rm meas}\nolimits}

   \def\no{\noindent}

   \def\ns{{ s}}

\theoremstyle{definition}
\newtheorem{df}{Definition}[section]
\theoremstyle{remark}
\newtheorem{rem}[df]{Remark}
\theoremstyle{plain}
\newtheorem{prop}[df]{Proposition}
\newtheorem{lemma}[df]{Lemma}
\newtheorem{teo}[df]{Theorem}

\newtheorem{cor}[df]{Corollary}

 \renewcommand{\theequation}{\thesection.\arabic{equation}}
 \newcommand{\sezione}[1]{\section{#1}\setcounter{equation}{0}}

\begin{document}


\title{Normalized positive solutions for Schr\"odinger equations with potentials in unbounded domains}

\author{
Sergio Lancelotti
\\
{\small\it Dipartimento di Scienze Matematiche, Politecnico di Torino}
\\
{\small\it Corso Duca degli Abruzzi n. 24, 10129 Torino, Italy}
\\
{\small\tt sergio.lancelotti@polito.it}
\\ \\
Riccardo Molle
\\
{\small\it Dipartimento di Matematica, Universit\`a di Roma ``Tor Vergata''}
\\
{\small\it Via della Ricerca Scientifica n. 1, 00133 Roma, Italy}
\\
{\small\tt molle@mat.uniroma2.it}
}

\maketitle


\begin{abstract}
The paper deals with the existence of positive solutions with prescribed $L^2$ norm for the Schr\"odinger equation
$$
-\Delta u+\lambda u+V(x)u=|u|^{p-2}u,\qquad u\in H^1_0(\Omega),\quad\int_\Omega u^2dx=\rho^2,\quad\lambda\in\R,
$$
where $\Omega=\R^N$ or $\R^N\setminus\Omega$ is a compact set, $\rho>0$, $V\ge 0$ (also $V\equiv 0$ is allowed), $p\in \left(2,2+\frac 4 N\right)$.
The existence of a positive solution $\bar u$ is proved when $V$ verifies a suitable decay assumption $(D_\rho)$, or if $\|V\|_{L^q}$ is small, for some $q\ge \frac N2$ ($q>1$ if $N=2$).
No smallness assumption on $V$ is required if the decay assumption $(D_\rho)$  is fulfilled.
There are no assumptions on the size of $\R^N\setminus\Omega$.
The solution $\bar u$ is a bound state and no ground state solution exists, up to the autonomous case $V\equiv 0$ and $\Omega=\R^N$.
\end{abstract}


{\em  \noindent \ \  A.M.S. Subject classification}: 35J50; 35J15; 35J60.




\bigskip

%

\sezione{Introduction and main results}


In this paper we study a class of problems of the form  
\beq
\label{P}
\tag{$P$}
\left\{
\begin{array}{ll}
-\Delta u+\lambda u+V(x)\, u=|u|^{p-2}u & \mbox{in $\Omega$}, \\
 \noalign{\medskip}
\lambda\in\R, \quad u\in H^1_0(\Omega) , \quad \displaystyle{\int_{\Omega} u^2\,dx=\rho^2},
\end{array}
\right.
\eeq
where $\Omega=\R^N$ or $\Omega$ is an exterior domain, that is $\R^N\setminus\Omega$ is a compact set, $N\ge 1$, $V$ is a given potential and $2<p<2+4/N$, namely the nonlinearity is superlinear and {\em mass subcritical}.
Here $\lambda$ will arise as a Lagrange multiplier related to the {\em mass constraint} $\|u\|_{L^2}=\rho>0$. 
We will focus on potentials that verify
\beq
\label{HV}
\begin{split}
    &V\in L^{q} (\R^N),\ q\in[\max(N/2, 1) ,\infty],\ \mbox{ with }q\neq 1\ \mbox{ if }N=2\\
    & V(x)\ge 0\ \mbox{ a.e. in }\R^N.
\end{split}
\eeq

\medskip

Problems of the form \eqref{P} arise from the  nonlinear Schr\"odinger and Klein-Gordon equations
\beq
\label{S}\tag{S}
i \Phi_t+\Delta\Phi+f(\Phi)=0
\eeq
\beq
\label{KG}\tag{KG}
\Phi_{tt}-\Delta\Phi+m\Phi=f(\Phi),
\eeq
where $\Phi:\R^N\times (0,+\infty)\to \C$.
If $f(\rho e^{i\theta})=f(\rho)e^{i\theta}$, $\rho,\theta\in\R$, then one can look for standing wave solutions of   \eqref{S} and \eqref{KG}, namely  solutions  of the form
\beq
\label{An}
\Phi(x,t)=e^{i\lambda t}u(x),\qquad x\in\R^N,\ t>0,
\eeq
where $u$ is a real function.
In the model case  $f(\Phi)=|\Phi|^{p-2}\Phi$ we consider, $u$ has to solve the equation in \eqref{P}.
We refer the reader to \cite{BeLi83ARMA1,Ca03book,CaLi82CMP} for more detailed physical motivations.

\medskip

If the frequency $\lambda$ in the ansatz \eqref{An} is fixed, and we consider the pure power model case with $p\in (2,2^*)$, where $2^*=\frac{2N}{(N-2)^+}$, then looking for solutions of \eqref{S} and \eqref{KG} corresponds to looking for critical points of the action functional
$$
E_\lambda (u)=\frac12\int_{\R^N}(|\D u|^2+[{\lambda}+V(x)]u^2)dx-\frac1p
\int_{\R^N}|u|^{p}dx,\qquad u\in H^1(\R^N).
$$
A very large number of works are devoted to this unconstrained problem, we only refer the reader to  \cite{AMbook06,C06MJ,St08book}  for a survey on almost classical results, and to the recent papers  \cite{DeSo20NA,21CoV} and references therein for new contributions.

\smallskip

Another point of view is to fix a priori the $L^2$ norm of the solution.
This point of view is particularly relevant because this quantity (the  {\em  mass} or the {\em charge} of the particle) is preserved along the time evolution.
In this case, the solutions $u$ correspond to the critical points of the energy functional
$$
E(u)=\frac{1}{2}\int_{\R^N} \left[|\nabla u|^2+V(x) u^2\right]\,dx -\frac{1}{p}\int_{\R^N} |u|^p\,dx\qquad u\in H^1(\Omega),
$$ 
constrained on 
$$
S_{\rho}:=\left\{u\in H^1(\R^N):\enskip |u|_2=\rho\right\},
$$
and the frequency $\lambda$ arises as a Lagrange multiplier.
Even if this fixed mass problem is classic (see for example \cite{CaLi82CMP}), only in the last decade it has been studied extensively and, in particular, very little has been done in the non-autonomous case (see \cite{AlJi22JGA,BMRV21,BaVa13AM,IkMi20CalcVar,Je97NA,JeLu19N,MRV22JDE,ZhZh22NODEA} and references therein).

\medskip

If $2+\frac 4N<p<2^*$, the so called mass-supercritical regime, then $E$ is not bounded from below on $S_\rho$, as follows by evaluating the functional over the fibers on $S_\rho$ of the type  $u_t:=t^{N/2} u_1(t\cdot)$, for fixed $ u_1\in S_\rho$ and $t>0$. 
So the problem cannot be addressed by minimization.
Jeanjean in the pioneering paper \cite{Je97NA} analyzed the autonomous problem, for more general mass-supercritical nonlinearities, and he recognized a mountain pass structure, related to the above-introduced fibers.  
The method developed in \cite{Je97NA} does not work in the non-autonomous case, even if the potential is radially symmetric.
The non-autonomous case in the mass-supercritical regime has been studied in \cite{BMRV21,MRV22JDE}, for not necessarily symmetrical potentials.
In \cite{BMRV21}  the authors found a mountain-pass solution if $V\ge 0$ is suitably small, while  \cite{MRV22JDE} concerns the case $V\le 0$ and the existence of two solutions is proved when $V$ is suitably small and the operator $-\Delta+V$  is not positive-definite.
  
\smallskip

In the mass-subcritical regime $2<p<2+\frac 4N$, the functional $E$ is bounded from below on $S_\rho$ by the Gagliardo-Nirenberg inequality (see \eqref{eGN}, \eqref{ciomp}). 

In \cite{IkMi20CalcVar}, Ikoma e Miyamoto  considered non-autonomous problems of the type \eqref{P} with more general mass-subcritical nonlinearities $u\mapsto f(u)$.
They assume  $V\le 0$ and $V(x)\to 0$ as $|x|\to\infty$ and prove, by concentration-compactness arguments, that a global minimum solution exists in $S_\rho$, for $\rho>\rho_0$, and does not exists if $0<\rho<\rho_0$, for a suitable $\rho_0\ge 0$.
Moreover, some sufficient conditions on $f$ and $V$ are provided to get $\rho_0=0$.
 
Alves and Ji in \cite{AlJi22JGA} consider some classes of potentials $V$ that do not vanish at infinity and where a global minimum for $E$ on $S_\rho$  again exists, for suitable $\rho$.
Namely,  in \cite{AlJi22JGA} the authors work on potentials that are $(V_1)$: 1-periodic in each variable, $(V_2)$: asymptotically 1-periodic, $(V_3)$: such that $\inf_{\R^N} V<\liminf_{|x|\to\infty}V(x)$, $(V_4)$: there exists $\mu_0>0$ such that $\meas\{V>\mu_0\}<\infty$ and $\interior(V^{-1}(0))\neq\emptyset$.

\medskip

This paper concerns the mass-subcritical case when $V\ge 0$ so no minimum solution exists, up to the autonomous case $V\equiv 0$ and $\Omega=\R^N$ (see Proposition \ref{prop:no_ground_state}). 
Moreover, we focus on domains that can be not only $\R^N$ but also exterior domains, answering the question raised in \cite{MRV22JDE} for mass-supercritical problems, whether the existence of bound state solution can be treated in exterior domains as in the whole space. 
To the best of our knowledge, these issues are not addressed in previous papers.
Concerning exterior domains, it is worth mentioning the recent paper \cite{ZhZh22NODEA} where the authors consider the autonomous problem $V\equiv 0$ in dimension $N\ge 3$ and they found the existence of a bound state solution if the size of ``hole'' $\R^N\setminus\Omega$ is small.

\smallskip

The main results are the following:

\begin{teo}
\label{T}
Let  $N\ge 2$,   $\Omega=\R^N$ or $\R^N\setminus\Omega$ compact,   $\rho>0$.
If  $V$ satisfies \eqref{HV} and
\beq
\label{HVas}
\tag{$D_\rho$}
\int_{\R^N} V(x)\,|x|^{N-1} e^{d_\rho|x|}\,dx<\infty,
\end{equation}
where
\beq
\label{703}
d_\rho=\left[2^{\left(1-\frac{p-2}{4-N(p-2)}\right)}\sqrt\lambda_1\right]\,  \rho^\frac{p-2}{4-N(p-2)},
\eeq
(see \eqref{Pkrho2} for the constant $\lambda_1$), then there exists a  solution $(\lambda,\bar u)$ of  \eqref{P}  such that $\lambda>0$ and $\bar u\ge 0$.

\end{teo}

\begin{teo}
\label{T2}
Let $N\ge 2$,   $\Omega=\R^N$ or $\R^N\setminus\Omega$ compact,   $\rho>0$.
There exists $L =L(q,\Omega,\rho)>0$ such that if $V$ satisfies \eqref{HV}, with
\beq
\label{infty}
  V(x)\lo 0\ \mbox{ as }|x|\to\infty\quad \mbox{if }V\in L^\infty(\Omega),
 \eeq
and  $\|V\|_{q}<L$, then   problem \eqref{P} has a solution $(\lambda,\bar u)$, verifying $\lambda>0$ and $\bar u\ge 0$.\end{teo}

The case $N = 1$ has its own specificity and will be dealt with in Section 6.

\medskip

A priori, a nonnegative solution $\bar u$  of \eqref{P} belongs to $H^1(\Omega)$, so we cannot say that $\bar u>0$, $\forall x\in\Omega$. 
Anyway, under mild assumptions, $\bar u$ turns out to be sufficiently regular to apply the Harnack inequality and therefore to be pointwise positive. 
In the following Proposition, we collect some regularity properties. 
Its proof is almost standard, so we only outline it in the Appendix.

\begin{prop}
\label{stab}
Let  $\bar u$ be a  solution of \eqref{P}. 
If $N\ge 2$, then
 
\begin{itemize}
\item[$ a)$]
if $V\in L^q_{\loc}(\Omega)$ for some $q>\frac N2$, then  $\bar u\in\cC^{0,\alpha}_{\loc}(\Omega)$, and if $\bar u$ is nonnegative then  $\bar u(x)>0$    $\forall x\in\Omega$;

\item[$ b)$] if $V\in L^q_{\loc}(\Omega)$ for some $q> N$, then $\bar u\in\cC^{1,\alpha}_{\loc}(\Omega)$; 

\item[$ c)$] if $V\in\cC^{0,\alpha}_{\loc}(\Omega)$, then $\bar u\in \cC^2(\Omega)$ and it is a classical solution.
\end{itemize}

If $N=1$,   then 
\begin{itemize}
\item[$ d)$]
if $V\in L^1_{\loc}(\Omega)$, then  $\bar u$ is continuously differentiable and if $\bar u$ is nonnegative then  $\bar u(x)>0$    $\forall x\in\Omega$;
\item[$ e)$] if $V\in L^q_{\loc}(\Omega)$ for some $q>1$, then $\bar u\in\cC^{1,\alpha}_{\loc}(\Omega)$.
\end{itemize}

\end{prop}

By Proposition \ref{stab},  the solutions given by Theorem \ref{T} and \ref{T2} are positive if assumption \eqref{HV} holds with $q>\frac N2$.

\medskip

The solution we find is a bound state solution.
Indeed, in Section 2 we verify that no ground state solution can exist:

\begin{prop} \label{prop:no_ground_state}
Assume that $V$ satisfies the assumptions of Theorem \ref{T} or of Theorem \ref{T2}.
If $V \not \equiv 0$ or $\Omega\neq\R^N$, then problem (\ref{P}) has no ground state solution.
\end{prop}

Actually, if $V$ satisfies both the assumptions of Theorem \ref{T} and of Theorem \ref{T2}, then the solutions provided by the Theorems are found exploiting the same topological configuration so it is reasonable to expect that they are the same solution.
Moreover, the topological characterization suggests that they have Morse index $N$.

\medskip

In  \eqref{HV}, the assumption that $V$ vanishes at infinity is necessary to get the compactness condition and cannot be dropped, by the following nonexistence result (see \cite{EsLi82PRSE}, \cite[Proposition 1.10]{MRV22JDE}, \cite[Theorem1.1]{CM10}).

\begin{prop}
\label{P:NE}
Let $p\in (2,2^*)$, $V\in   L^\infty(\R^N)$  and assume that there exists  $\frac{\partial V}{\partial\nu}\in L^{\ell}(\R^N)$ for some $\nu\in\R^N\setminus\{0\}$ and $\ell\in \left[\max(1,\frac{N}{2}) ,+\infty\right]$, $\ell\neq 1$ if $N=2$. 
If $\frac{\partial V}{\partial\nu}\ge 0$ and $\frac{\partial V}{\partial\nu}\not\equiv 0$, then problem
 \beq
 \label{Rs}
 -\Delta u+\lambda u-V(x)u=|u|^{p-2}u\qquad u\in S_\rho, \quad \lambda\in\R
\eeq
 has no solutions in $\cC^1(\R^N)\cap W^{2,2}(\R^N)$
\end{prop}
(see also Remark \ref{R:EL} and Proposition \ref{P:NEN=1} for the 1-dimensional case).

\medskip

Some remarks are in order, concerning the decay assumption \eqref{HVas} in Theorem \ref{T}.
It does not require any smallness assumption on $V$ and moreover  does not depend on the domain, hence if it is verified then problem \eqref{P} has a solution for every exterior domain $\Omega$.

On the other hand, \eqref{HVas} depends on $\rho$ by \eqref{703}, where $d_\rho\to\infty$ as $\rho\to\infty$ because $p<2+\frac 4N$.
As a consequence, a potential that verifies \eqref{HVas} for every $\rho>0$ has to decay faster than any exponential.
 
\smallskip
 
\begin{rem}
\label{aas}
By the  proof of Theorem \ref{T2} we verify the existence of a constant $\bar L=\bar L(N)>0$ such that 
\beq
\label{eaas0}
L(N/2,\R^N,\rho)=\bar L \qquad\forall \rho>0 .
\eeq
Instead, for $q>\frac N2$,
\beq
\label{eaas}
\lim_{\rho\to 0}L(q,\R^N,\rho)=0, \qquad \lim_{\rho\to\infty}L(q,\R^N,\rho)=\infty.
\eeq
As a consequence, if $\|V\|_{N/2}$ is suitably small, then problem \eqref{P} in $\R^N$ has a solution for every $\rho>0$, while if $V\in L^q(\R^N)$, for $q>N/2$, then there exists $\bar \rho=\bar \rho({\|V\|_q})\ge 0$ such that  problem \eqref{P} has a solution for every $\rho>\bar\rho$ (see \eqref{Sergio}).
\end{rem}
In Theorems 1.1. and 1.2 of \cite{MRV22JDE} a mountain pass solution is found in the mass supercritical case with negative potential, under smallness assumptions on $V$ similar to the ones considered here. 
In that case, if $q=N/2$ there is a uniform bound as in \eqref{eaas0} while for $q>N/2$ the limits in \eqref{eaas} are reversed. 

\medskip

If $\rho$ is fixed, we will observe that  
$$
\lim_{r(\Omega)\to\infty} L(q,\Omega,\rho)=0,
$$
where
\beq
\label{7}
r(\Omega)=\sup\{r\in\R^+\ :\ B_r(y)\subset \R^N\setminus\Omega\mbox{ for some } y\in\R^N\}
\eeq
(see Remark \ref{ras}).
 Hence, there is no potential $V\not\equiv 0$ such that Theorem \ref{T2} applies for every exterior domain $\Omega$, for  $\rho>0$ fixed.

\smallskip

If we consider $V\equiv 0$, then both Theorems \ref{T} and \ref{T2} apply and provide the following result.

\begin{cor}
\label{Caut}
Let $N\ge 2$,  $\R^N\setminus \Omega$ compact and $\rho>0$, then there exist $\lambda>0$ and $\bar u\in   \cC^2(\Omega)$, $\bar u(x)>0$, $\forall x\in\Omega$, such that 
$$
\left\{
\begin{array}{ll}
-\Delta \bar u+\lambda \bar u =\bar u^{p-1} & \mbox{in $\Omega$},  
 \\
 \noalign{\medskip}
\lambda\in\R, \quad \bar u\in H^1_0(\Omega),\quad \displaystyle{\int_{\Omega} \bar u^2\,dx=\rho^2}.
\end{array}
\right.
$$
\end{cor}

Corollary \ref{Caut} extends the result of \cite{ZhZh22NODEA} to the dimension 2 and to every exterior domains. 
The proof in \cite{ZhZh22NODEA} cannot be extended to this more general framework because it does not work in dimension 2 and when the size of $\R^N\setminus\Omega$ is large.

\medskip

We prove Theorems \ref{T} and \ref{T2} by variational methods, looking for bound state solutions.
The analysis of the compactness presents a lot of difficulties related to the unboundedness of the domain, that is not assumed to be symmetric.
In order to recover a local compactness condition, we first see that the Lagrange multiplier $\lambda$ related to a Palais-Smale sequence in a negative range is positive, and then we employ a splitting Lemma from \cite{BC87ARMA} for PS sequences of $E_\lambda$.
With this decomposition in hands, we perform in Proposition \ref{PPS} a sharp fine estimate of the first energy interval $I$ above the infimum where the compactness condition holds (see Remark \ref{MC}).

\smallskip

The topological argument relies on min-max techniques that make use of a barycentric map and the Brouwer
degree. 
Since we have no smallness assumption on $\R^N\setminus\Omega$, and in Theorem \ref{T} we have no bound on any Lebesgue norm of the potential, a major difficulty is to work in the compactness interval $I$, in the min-max procedure.
To overcome this problem, we will proceed by analyzing the energy interaction of positive solutions $w_1$ and $w_2$ of some suitable  ``problems at infinity'', such that $\|w_1\|_{L^2}^2+\|w_2\|_{L^2}^2=\rho^2$.
This idea is inspired by the unconstrained case (see \cite{CP95,LM20}), where the problem at infinity is univocally determined by the choice of $\lambda$.
Here, the need to choose different functions $w_1,w_2$ prevents the use of the arguments developed in \cite{CP95,LM20}) and requires the introduction of different and more refined estimates.

\smallskip

To verify that the solutions we find do not change sign, we prove in Proposition \ref{Psegno} that the energy of every solution that changes sign is not in the energy interval we are working in.   
In particular, Proposition \ref{Psegno} and Corollary \ref{Psgnmassa} give information also on changing sign solutions of the autonomous problem in $\R^N$.
We point out that to get nonnegative solutions here we could exploit the symmetry of the functional and work near the cone of the positive functions, by using \cite[Theorem 4.5]{Gh93book} and proceeding as in \cite{BMRV21}.
The advantage of this other approach is a simplification of the proof of the compactness condition, because the lack of compactness in such a case comes only from the positive solution of the problem at infinity.
On the other hand, the approach we employ here allows us both to bound from below the energy of the changing sign solutions and to gain a more general analysis of the Palais-Smale sequences.

\medskip

The paper is organized as follows: in Section 2 we introduce some preliminary results, in Section 3 we prove the local compactness condition and  Section 4 is devoted to the sharp energy estimates that are necessary in Section 5 to prove Theorems \ref{T} and \ref{T2}, Section 6 concerns the case $N=1$ and in the Appendix we give a sketch of the proof of the regularity proposition, with some references for detailed proofs.
 
\pagebreak 
 

\sezione{Notations, variational framework and preliminary results}


Throughout the paper we make use of the following notation:

{\small
 \begin{itemize}

\item $2_c:=2+\frac 4 N$,\quad  $\ns=\frac2N\frac{p-2}{2_c-p}$.
 
\item $L^q(\mathcal{O})$, $1\leq q \leq \infty$, $\mathcal{O}\subseteq
  \R^N$ a measurable set,
  denotes the Lebesgue space, the norm in $L^q(\mathcal{O})$ is denoted
  by $|\cdot|_q$ if there is no ambiguity.

\item For $u\in H^1_0(\Omega)$ we denote by $u$ also the function in
  $H^1(\R^N)$ obtained setting $u\equiv 0$ in $\R^N\setminus\Omega$.
  
\item For any $R>0$ and for any $z\in \R^N$, $B_R(z)$
  denotes the closed ball of radius $R$ centered at $z,$ and for any
  measurable set $ \mathcal{O} \subset \R^N, \ |\mathcal{O}|$ denotes
  its Lebesgue measure.

\item $H^{1}(\R^{N})$ is the usual Sobolev space endowed with
  the standard  norm
\begin{displaymath}
 \|u\|:=\left[\int_{\R^N}\left(|\nabla u|^{2}+ u^{2}\right)dx\right]^{\frac 12}.
\end{displaymath}

\item $c,c', C, C', C_i,\ldots$ denote various positive constants that can also vary from one line to another.
\item $o(f)$ and $O(f)$ denote the usual Landau notations: $\frac {o(f)}{f}\to 0$ as $f\to 0$ and $|O(f)|\le C  |f|$ for some positive constant $C$.

\end{itemize}
}

We will find solutions $\bar u$ of problem $(P)$ as critical points of the functional $E$ constrained on $S_{\rho}$. 
If $\lambda$ is the Lagrange multiplier related to $\bar u$, then $\bar u$ is also a free critical point of the related free functional $E_\lambda$

\medskip

Let us assume $\lambda>0$ and recall some well-known properties of the limit problem, for $V\equiv 0$ and $\Omega=\R^N$,
\beq
\label{Pinfty}
\tag{$P_\infty$}
\left\{
\begin{array}{ll}
-\Delta u+\lambda_{\infty} u=|u|^{p-2}u & \mbox{in $\R^N$}, \\
\noalign{\medskip}
\lambda_{\infty} \in\R,
\quad
u\in S_\rho.
\end{array}
\right.
\eeq
\eqref{Pinfty} has a unique positive solution $w\in H^1(\R^N)$, up to translations, which is radial, radially decreasing, and belongs to $C^2(\R^N)$.
The function $w$ verifies the minimality property
\beq
\label{em}
m:=E_\infty(w)=\min_{u\in S_\rho} E_\infty(u),
\eeq
where
$$
E_\infty(u)=\frac{1}{2}\int_{\R^N} |\nabla u|^2\,dx -\frac{1}{p}\int_{\R^N} |u|^p\,dx.
$$

Correspondingly, the solutions of \eqref{Pinfty} are also free critical points of the limit functional
$$
E_{\lambda,\infty}(u)=E_\infty (u)+\frac{\lambda}{2}\int_{\R^N} u^2dx\qquad u\in H^1(\R^N).
$$

Moreover,
\beq
\label{<0}
m<0\quad\mbox{ and }\quad \lambda_\infty>0
\eeq
and there exists $c_1>0$ such that
\begin{equation} \label{eqn:decrescita_w}
w(|x|)\, e^{\sqrt{\lambda_{\infty}}|x|}\,|x|^{\frac{N-1}{2}} \to c_1 \quad \mbox{as $|x|\to \infty$}, \\
\end{equation}
\begin{equation} \label{eqn:decrescita_w'}
w'(|x|)\, e^{\sqrt{\lambda_{\infty}}|x|}\,|x|^{\frac{N-1}{2}} \to -c_1\sqrt{\lambda_{\infty}} \quad \mbox{as $|x|\to \infty$}. \\
\end{equation}

Inequality $m<0$ in \eqref{<0} follows choosing $u_1\in S_\rho$ and testing \eqref{em} on $u_t:=t^{N/2}u_1(t\cdot)\in S_\rho$, $t>0$, taking into account $p<2_c$.
The positivity of  $\lambda_\infty$ comes from Pohozaev and Nehari identities (see also Proposition \ref{Psegnolambda}).
For the properties of $w$ we refer the reader to \cite{BaLi97AIHP,St77CMP,BeLi83ARMA1,Kw89ARMA,GNN81}.

For any $k>0$, let us denote by $w_{k\rho^2}$ the positive solution of
\beq
\label{Pkrho2}
\left\{
\begin{array}{ll}
-\Delta u+\lambda_{k\rho^2}\, u=|u|^{p-2}u & \mbox{in $\R^N$}, \\
\noalign{\medskip}
\lambda_{k\rho^2}\in\R,\quad u\in S_{ \sqrt k\, \rho},  \\
\end{array}
\right.
\eeq
where $w_{k\rho^2}$ verifies
\beq
\label{699}
m_{k\rho^2}:=E_{\infty}(w_{k\rho^2})=\inf_{u\in S_{\sqrt k\, \rho}} E_{\infty}(u).
\eeq
It turns out that
\beq
\label{700}
w_{k\rho^2}(x)=k^{\frac{\ns}{p-2}}w\left(k^{\frac{\ns}{2}}x\right),
\quad
m_{k\rho^2}=k^{1+\ns}\,m<0,
\quad
\lambda_{k\rho^2}=k^{\ns}\,\lambda_{\infty}>0
\eeq
where 
$$
\ns:=\frac2N\frac{p-2}{2_c-p}
$$ 
(see notations). 
From $p<2_c$ it follows  $s>0$.
By \eqref{700}, in particular we have
\beq
\label{702}
\lambda_{\rho^2}=\lambda_\infty =\rho^{s}\lambda_1.
\eeq
 For $k=0$ we set $w_0=0$ and for $k=1$ we simply write $w_{\rho^2}=w$.

Moreover, for any $k>0$, by (\ref{eqn:decrescita_w}), (\ref{eqn:decrescita_w'}) and \eqref{700} we have
\begin{equation} \label{eqn:decrescita_w_k}
w_{k\rho^2}(|x|)\, e^{\sqrt{k^{\ns}\lambda_{\infty}}|x|}\,|x|^{\frac{N-1}{2}} \to c_k \quad \mbox{as $|x|\to +\infty$}, \\
\end{equation}
\begin{equation} \label{eqn:decrescita_w_k'}
w_{k\rho^2}'(|x|)\, e^{\sqrt{k^{\ns}\lambda_{\infty}}|x|}\,|x|^{\frac{N-1}{2}} \to -c_k\,\sqrt{k^{\ns}\,\lambda_{\infty}} \quad \mbox{as $|x|\to +\infty$}\\
\end{equation}
where
$$
c_k=c_1\, k^{\ns\left(\frac {1}{p-2}-\frac{N-1}{4}\right)}.
$$
From $p<2_c$ there follows $\frac {1}{p-2}-\frac{N-1}{4}>0$, so that
$$
c_k\to 0\quad\mbox{ as } k\to 0,
\qquad
c_k\to c_1\quad\mbox{ as } k\to 1.
 $$

{\noindent{\em Proof of Proposition \ref{prop:no_ground_state}.}}\hspace{2mm}
Let us prove that
\beq
\label{1229}
\inf_{S_{\rho}} E=m
\eeq
and the infimum is not achieved.

Obviously we have that ${\inf_{S_{\rho}} E\geq m}$.
To show that the equality holds, let us consider the sequence $(u_n)_n$ defined by $u_n=w(x-y_n)$, where $(y_n)_n$ is a sequence in $\R^N$ such that $|y_n|\to +\infty$.
By (\ref{HV}) and (\ref{eqn:decrescita_w}), and taking into account \eqref{HVas} or \eqref{infty}, we have that
$$
\lim_{n\to +\infty} \int_{\R^N} V(x) u_n^2(x)\,dx =0,
$$
which implies $\displaystyle{\lim_{n\to +\infty} E(u_n)=m}$.

Now,  assume by contradiction that the infimum ${\inf_{S_{\rho}} E}$ is achieved by a function $\bar u$.
Then, by
\beq
\label{cov}
m\le E_\infty(\bar u)\le E(\bar u)=m,\qquad \bar u\in S_\rho,
\eeq
and by the uniqueness of the minimizers of \eqref{em},  we should have $\bar{u}=\pm w(x-\bar{y})$ for a suitable $\bar y\in\R^N$.
Since  $w(x-\bar{y})>0$ for all $x\in\R^N$, we can deduce that $\Omega$ must be the entire space $\R^N$ and by \eqref{cov}
we get
$$
\int_{\R^N} V(x) \bar{u}^2(x)\,dx =0,
$$
that implies $V \equiv 0$ because $V\ge 0$.
So we are in contradiction with $\Omega\neq\R^N$ or $V\not\equiv 0$.

\qed

Next proposition states the positivity of the Lagrange multipliers in \eqref{P}.

\begin{prop}
\label{Psegnolambda}
Assume that  $u\in S_\rho$ and $\lambda\in \R$ solve \eqref{P}.
If $E(u)<0$, then $\lambda>0$.
\end{prop}

\proof
We have
$$
\int_\Omega \D u\cdot \D v\, dx +\lambda\int_\Omega u v\, dx +\int_\Omega V(x)u v\, dx-\int_\Omega |u|^{p-2}uv\, dx=0\qquad \forall v\in H^1_0(\Omega).
$$

Then,
\beq
\label{1813}
\lambda\rho^2=|u|_p^p-|\D u|_2^2-\int_\Omega V(x)u^2  dx =-p E(u)+\frac{p-2}{2}|\D  u|_2^2+\frac{p-2}{2}\int_\Omega V(x)u^2  dx>0.
\eeq

\qed

Before concluding this section, let us find out some features of the changing sign solutions.

 \begin{prop}
 \label{Psegno}
 Let $u\in S_\rho$ be a solution of \eqref{P} ($(P_\infty)$).
 If $u$ changes sign, then $E(u) > 2^{-s}m$ ($E_\infty(u) > 2^{-s}m$).
  \end{prop}

 \proof
Since $m<0$, we can assume $E(u)<0$.
According to Proposition \ref{Psegnolambda}, let $\lambda>0$ be the Lagrange multiplier corresponding to the solution $u$ and let $w_\infty\in H^1(\R^N)$ be the positive solution of $-\Delta v+\lambda v=v^{p-1}$.
By \eqref{700},
\beq
\label{1104}
E_\infty(w_\infty)=\min_{S_{|w_\infty|_2}}E_\infty=\left(\frac{|w_\infty|_2^2}{\rho^2}\right)^{1+s}m.
\eeq

Let us write $u=u^+-u^-$, where $u^\pm=\max\{\pm u,0\}$, $\rho^2=|u^+|_2^2+|u^-|_2^2$.
Since $u^+$ and $u^-$ are on the Nehari manifold $\cN$ corresponding to $E_\lambda$, then $(E_\lambda)_{|_\cN}(u^\pm) > E_{\lambda,\infty}(w_\infty)$ (see \cite{BC87ARMA} and \cite[Proposition 2.1]{MMP00PRSE}).
Namely, by \eqref{1104},
$$
E(u^\pm)> \left(\frac{|w_\infty|_2^2}{\rho^2}\right)^{1+s}m+\frac\lambda 2(|w_\infty|_2^2-|u^\pm|_2^2).
$$

Hence, the claim read as
$$
E(u)=E(u^+)+E(u^-)>2\, \left(\frac{|w_\infty|_2^2}{\rho^2}\right)^{1+s}m+\frac\lambda 2(2|w_\infty|_2^2-\rho^2)\ge 2^{-s}m,
$$
and, equivalently, as
\beq
\label{1333}
\left[\left(\frac{2|w_\infty|_2^2}{\rho^2}\right)^{1+s}-1\right]2^{-s} m
+\frac\lambda 2\rho^2\,\left[\left(\frac{2|w_\infty|_2^2}{\rho^2}\right)-1\right]\ge 0.
\eeq
Since $(\lambda,w_\infty)$ solves  \eqref{Pkrho2} for $k=\frac{|w_\infty|^2_2}{\rho^2}$, the following system is fulfilled
\beq
\label{sys}
\left\{
\begin{array}{lcl}
\vspace{2mm}
\frac 12 |\D w_\infty|_2^2-\frac 1p|w_\infty|_p^p=m_{|w_\infty|_2^2} &&\mbox{\em Energy of $w_\infty$}\\
\vspace{2mm}
|\D w_\infty|_2^2+\lambda | w_\infty|_2^2-|w_\infty|_p^p=0&&\mbox{\em Nehari}\\
\vspace{2mm}
\frac{N-2}{2}|\D w_\infty|_2^2+\frac N 2\lambda|w_\infty|_2^2-\frac{N}{p}|w_\infty|_p^p=0&&\mbox{\em Pohozaev.}
\end{array}
\right.
\eeq

Solving \eqref{sys} and taking into account \eqref{700} we get
\beq
\label{1334}
\frac \lambda 2\,\rho^2
=-\frac{2p-N(p-2)}{4-N(p-2)}\, \left(\frac {|w_\infty|_2^2}{\rho^2}\right)^{s}m.
\eeq

Putting \eqref{1334} in \eqref{1333}, we have to verify
$$
\left[\left(\frac{2|w_\infty|_2^2}{\rho^2}\right)^{1+s}-1\right]
-\frac{2p-N(p-2)}{4-N(p-2)} \,\left(\frac {2|w_\infty|_2^2}{\rho^2}\right)^s\left[\left(\frac{2|w_\infty|_2^2}{\rho^2}\right)-1\right]\le 0.
$$

Then the assertion is proved because
$
\frac{2p-N(p-2)}{4-N(p-2)}=1+s
$
and
$$
[t^{1+s}-1]-(1+s) \, t^s[t-1]\le 0\qquad\forall t\ge 0.
$$

\qed

\begin{cor}
\label{Psgnmassa}
If $u$ is a changing sign solution of
\beq
\label{1044}
-\Delta v+\lambda v=|v|^{p-2}v\qquad v\in H^1(\R^N),
\eeq
then $|u|_2^2> 2|w_\infty|_2^2$, where $w_\infty$ denotes the positive solution of \eqref{1044}.
\end{cor}

Clearly, here  $\lambda>0$ by Pohozaev identity.

\proof
Let us call $E(u)=M$.
From system \eqref{sys} applied to $u$ and $w_\infty$ it follows
\beq
\label{1111}
|u|_2^2=2\, \frac{1+s}{\lambda}\, (-M),\qquad |w_\infty|_2^2=2\, \frac{1+s}{\lambda}\, \left(-m_{|w_\infty|_2^2}\right).
\eeq
By Proposition \ref{Psegno}, $M>2^{-s}m_{|u|_2^2}$, then, by \eqref{700} and \eqref{1111},
$$
|u|_2^2 < 2\, \frac{1+s}{\lambda}\, 2^{-s} \left(-m_{|u|_2^2}\right)= 2\, \frac{1+s}{\lambda}\, 2^{-s}\left(\frac{|u|_2^2}{|w_\infty|_2^2}\right)^{1+s}  \left(-m_{|w_\infty|_2^2}\right)=2^{-s} |w_\infty|_2^2 \left(\frac{|u|_2^2}{|w_\infty|_2^2}\right)^{1+s},
$$
that completes the proof.

\qed


\sezione{The compactness condition}


\begin{prop}
\label{PPS}
Let $(u_n)_n$ be a Palais-Smale sequence at the level $c$ for $E$ constrained on $S_\rho$.
If $c\in(m,2^{-s}m)$ then there exists a critical point $u_0\in S_\rho$ such that $u_n\to u_0$, as $n\to\infty$.
\end{prop}

To prove Proposition \ref{PPS}, we will use the well known splitting lemma of Benci and Cerami for the unconstrained problem
(\cite[Lemma 3.1]{BC87ARMA}).

\begin{lemma}
\label{LBC}
Let $\lambda>0$ and let $(u_n)_n$ in $H^1(\R^N)$ be a Palais-Smale sequence for $E_\lambda$.
Then there exist a critical point $u_0$ of $E_\lambda$, an integer $h\ge 0$, $h$ non-trivial solutions $w^1,\dots,w^h\in H^1(\R^N)$ to the limit equation
$$-\Delta v+\lambda v=|v|^{p-2}v$$
and $h$ sequences $(y_n^j)_n\subset\R^N$, $1\le j\le h$, such that $|y_n^j|\to\infty$ as $n\to\infty$, and
\begin{equation}
\label{ss}
u_n=u_0+\sum_{j=1}^h w^j(\cdotp-y^j_n)+o(1)\qquad\text{strongly in $H^1(\R^N)$,}
\end{equation}
up to a subsequence.
Moreover, we have
 \begin{equation}
\label{se}
E_\lambda (u_n)\to E_\lambda (u_0)+\sum_{j=1}^h E_{ \lambda,\infty}(w^j)\qquad\text{as $n\to\infty$,}
\end{equation}
and
 \begin{equation}
  \label{sn}
|u_n|_2^2=|u_0|_2^2+\sum_{j=1}^h |w^j|_2^2+o(1).
\end{equation}
\end{lemma}

{\em Proof of Proposition \ref{PPS}}.\
We claim that $(u_n)_n$ is bounded in $H^1_0(\Omega)$.
Indeed, $u_n\in S_\rho$, $\forall n\in \N$, and by the Gagliardo-Nirenberg inequality 
\begin{equation}
\label{eGN}
|u|_p^p\leq C_{GN} |u|_2^{p-\frac{N(p-2)}{2}}|\D u|_2^{\frac{N(p-2)}{2}}.
\end{equation}
(see \cite[Theorem 12.83]{LeoniBook}) we have  
\beq
\label{ciomp}
c+o(1)=E(u_n)\ge\frac 12 |\D u_n|^2_2-C |\D u_n|_2^{\frac{N(p-2)}{2}},
\eeq
with $\frac{N(p-2)}{2}<2$ because $p<2_c$.

Since $(u_n)_n$ is a constrained PS-sequence, there exists a sequence $(\lambda_n)_n$ in $\R$ such that
\beq
\label{1812}
\int_\Omega \D u_n\cdot \D v\, dx +\lambda_n\int_\Omega u_n v\, dx+ \int_\Omega V(x)\, u_n v\, dx -\int_\Omega |u_n|^{p-2}u_nv\, dx=o(1)\|v\| \qquad \forall v\in H^1_0(\Omega).
\eeq
Setting $v=u_n$ in \eqref{1812}, and taking into account that $(u_n)_n$ is bounded in $H^1$, we can argue as in \eqref{1813} obtaining
\beq
\label{1821}
\lambda_n\rho^2
=-pE(u_n)+\frac{p-2}{2}|\D u_n|_2^2 +\frac{p-2}{2}\int_\Omega V(x)\, u_n^2 dx  +o(1) \ge -p\, c+o(1)>0,
\eeq
for large $n$.
Since $(u_n)_n$ is bounded in $H^1_0(\Omega)$, the first relation in  \eqref{1821} implies that  the sequence $(\lambda_n)_n$ is bounded. 
Moreover, from  \eqref{1821}  we infer also that $\lambda_n\ge c>0$, for a suitable constant $c$.
Then we can assume that $\lambda_n\to\lambda>0$.
Hence, by \eqref{1812} we are in position to apply Lemma \ref{LBC} and we can decompose $u_n$ according to \eqref{ss}.

\medskip

$(I)$\ If $h=0$, we are done.

\medskip

So, we assume by contradiction $h\ge 1$ and we are going to show that then
 \begin{equation}
\label{711}
E  (u_0)+\sum_{j=1}^h E_{\infty}(w^j)\ge 2^{-s}m,
\end{equation}
up to the case $u_0=0$, $h=1$ and $w^1>0$, that arises for
\beq
\label{eT}
\lim_{n\to\infty}E(u_n)=m.
\eeq
Once \eqref{711} is proved, the contradiction comes out, because  \eqref{se} and \eqref{sn} provide $c=E  (u_0)+\sum_{j=1}^h E_{\infty}(w^j)$, and $c\in(m,2^{-s}m)$ by assumption.

\medskip

$(II)$\  If $u_0\equiv 0$ and $h=1$  occours, then   $|w^1|_2=\rho$, by \eqref{sn}.
Hence,  if $w^1>0$ then  $c=m$ by \eqref{se} and \eqref{em}, so \eqref{eT} is proved.
We observe that this is the only case in which \eqref{eT} holds, up to the autonomous case $\Omega=\R^N$ and $V\equiv 0$, when also the case $u_0>0$ and $h=1$ verifies \eqref{eT}.
On the other hand, if  $w^1$ is a changing sign solution, then  $c>2^{-s} m$ by Proposition \ref{Psegno}.

\medskip

$(III)$\  If $u_0\equiv 0$ and $h\ge 2$, we proceed by induction.
For $h=2$, we get \eqref{711} arguing exactly as in the proof of Proposition \ref{Psegno}, with $w^1$ in place of $u^+$ and $w^2$ in place of $u^-$.
Observe that the case $w^1=w^2>0$ is the only case when the equality holds in \eqref{711}.

  Suppose now $h\ge 3$ and \eqref{711} holds for $h-1$, namely for every $\rho_1>0$
\beq
\label{1408}
\sum_{j=1}^{h-1} E_{\infty}(w^j)\ge 2^{-s}m_{\rho_1^2},
\eeq
whenever $w^j\not\equiv 0$, for all $j\in\{1,\ldots,h-1\}$, and  $\sum_{j=1}^{h-1}|w^j|_2^2=\rho_1^2$.
Then, let us prove \eqref{711} for $h$.
We can assume that $ |w^h|_2^2\le \rho^2/3$.
Taking into account \eqref{1408}  and \eqref{700}, we have
\begin{equation*}
\begin{split}
\sum_{j=1}^hE_\infty(w^j)& \ge
2^{-s} m_{(\rho^2-|w^h|_2^2)} +E_\infty(w^h)\\
&
\ge
2^{-s} m_{(\rho^2-|w^h|_2^2)} +m_{|w^h|_2^2} \\
&=
2^{-s}\left(\frac{\rho^2-|w^h|_2^2}{\rho^2}\right)^{1+s}m+\left(\frac{|w^h|_2^2}{\rho^2}\right)^{1+s}m.
\end{split}
\end{equation*}
Hence, it is sufficient to verify that
\beq
\label{1429}
2^{-s}\left(\frac{\rho^2-|w^h|_2^2}{\rho^2}\right)^{1+s}m+\left(\frac{|w^h|_2^2}{\rho^2}\right)^{1+s}m\ge 2^{-s}m.
\eeq
Inequality \eqref{1429} is equivalent to
\beq
\label{1637}
2^s t^{1+s}+(1-t)^{1+s}\le 1,
\eeq
where $t:=\frac{|w^h|_2^2}{\rho^2}\in(0,1/3]$.
Since inequality \eqref{1637} holds for every $t\in[0,1/3]$, estimate \eqref{711} is proved for $u_0\equiv 0$.

\medskip

$(IV)$\ If $u_0\not\equiv 0$, we can proceed as in the previous steps: first by considering the case $h=1$ and arguing as in the proof of Proposition \ref{Psegno}, and then finishing the proof by induction.

 \qed

\begin{cor}
\label{Cmin}
If $(u_n)_n$ in $S_\rho$ satisfies $\lim\limits_{n\to\infty}E(u_n)=m$, then there exists a sequence $(y_n)_n$ in $\R^N$ such that
\beq
\label{emm}
u_n(x)=w(x-y_n)+o(1)\qquad\mbox{ in }H^1.
\eeq
If $\Omega\neq\R^N$ or $V\not\equiv 0$, then $|y_n|\lo\infty$, as $n\to\infty$.
\end{cor}

\proof
By the Ekeland variational principle, there exists a PS sequence $(v_n)_n$ for $E$ constrained on $S_\rho$ such that  $\lim\limits_{n\to\infty}E(v_n)=m$ and $u_n=v_n+o(1)$ in $H^1$ (see  (\cite[Proposition 5.1]{Ek74JMAA}), then we can assume that $(u_n)_n$ is PS sequence.
Then \eqref{emm} is a direct consequence of \eqref{eT} and part $(II)$ in the proof of Proposition \ref{PPS}.

\smallskip

If $y_n\to\bar y\in\R^N$, up to a subsequence, then $u_n\to w(\cdot-\bar y)$ in $H^1$ and a.e., so that $\Omega=\R^N$ because $w(x)>0$, $\forall x \in \R^N$.
Moreover
\begin{equation*}
\begin{split}
m=\lim_{n\to\infty}E(u_n)&=\frac12\int_{\R^N} |\D w(x-\bar y)|^2dx+\frac12\int_{\R^N} V(x)\,  w(x-\bar y)^2dx-\frac 1p\int_{\R^N} w^p(x-\bar y)\, dx\\
&=m+\frac12\int_{\R^N} V(x)\,  w(x-\bar y)^2dx
\end{split}
\end{equation*}
implies $V\equiv 0$, again because $w(x)>0$, $\forall x\in\R^N$.

\smallskip

So the proof is completed.

\qed

\begin{rem}
\label{MC}
By $(III)$ in the proof of Proposition \ref{PPS} we see that, for every sequence $(y_n)_n$ in $\R^N$ such that $\lim_{n\to\infty}|y_n|=\infty$, 
$$
E\left(w_{ {\rho^2}/{2}}(\cdot-y_n)+w_{ {\rho^2}/{2}}(\cdot+y_n)\right)=2^{-s}m.
$$
Then the sequence $\big(w_{{\rho^2}/{2}}(\cdot-y_n)+w_{{\rho^2}/{2}}(\cdot+y_n)\big)_n$ turns out to be a not relatively compact PS-sequence at the level $2^{-s}m$, showing that the compactness interval $\big(m,2^{-s}m\big)$ cannot be extended.
\end{rem}


\sezione{Energy estimates}


\subsection{An upper bound}

If $\R^N\setminus\Omega$ is a non empty compact set contained in $B_{R-1}(0)$, let us introduce the cut-off function $\vartheta$,  verifying
\beq
\label{etheta}
\vartheta\in
\cC^\infty(\R^N,[0,1]),\qquad
\left\{\begin{array}{ll}
\vartheta(x)=1&\mbox{ if } |x|\ge R\\
\vartheta(x)=0&\mbox{ if }x\in\R^N\setminus\Omega.
\end{array}\right.
\eeq
If $\Omega=\R^N$, we agree that  $\vartheta\equiv 1$ on $\R^N$.
Let us set $\Sigma=\partial B_2(e_1)$, where $e_1$ is the first element
of the canonical basis of $\R^N$, and for any $r>0$ define the map
$\psi_r: [0,1]\times\Sigma\longrightarrow H^1_0(\Omega)$ by
\beq
\label{100}
\psi_{r}[t,z](x)= \rho\,\frac{\vartheta(x)\left[w_{t\rho^2}(x-r z)+w_{(1-t)\rho^2}(x-r e_1)\right]}
{\left|\vartheta(\cdot)\left[w_{t\rho^2}(\cdot-r z)+w_{(1-t)\rho^2}(\cdot-r e_1)\right]\right|_2}.
\eeq

\begin{prop} \label{lem:A_varepsilon<2m}
Suppose that $V$ verfies \eqref{HV} and \eqref{HVas}, then
\begin{itemize}
\item[$(a)$]
 there exist $ \overline{r}>0$ such that for any
$r>\overline{r}$
\beq
\label{A_r}
{\cal A}_r=\max\left\{E\left(\psi_{r}[t,z]\right):\enskip t\in [0,1],\,\, z\in\Sigma\right\} <2^{-\ns}\,m;
\eeq
\item[$(b)$]
for every $\e>0$ there exists $r_\e>0$ such that for any $r>r_\e$
\beq
\label{936}
{\cal L}_r=\max\left\{E\left(\psi_{r}[1,z]\right)\ :\   z\in\Sigma\right\} \leq m+\e.
\eeq
\end{itemize}
\end{prop}

\begin{rem}
\label{remV=0}
In Proposition \ref{lem:A_varepsilon<2m}, $V\equiv 0$ is allowed.
\end{rem}

Before proving Proposition \ref{lem:A_varepsilon<2m}, let us recall two
technical lemmas. For the proof of Lemma \ref{L4.4} we refer to
\cite{CP95} while the proof of Lemma \ref{BL} is in  \cite{BaLi90RMI} (see
also Lemma 2.9 in \cite{CM16}).

\begin{lemma}
\label{L4.4}
For all $a,b\ge 0$, for all $p\ge 2$, the following relation holds true
\[
(a+b)^{p}\ge a^{p}+b^{p}+(p-1)(a^{p-1}b+ab^{p-1}).
\]
\end{lemma}

\begin{lemma}
\label{BL}
If $g \in L^\infty (\R^N)$ and $h\in L^1(\R^N)$ are
such that, for some $\alpha\ge 0$, $b\ge 0$, $\gamma\in \R$
\beq
\label{eBLg}
\lim_{|x|\to \infty} g(x)\,e^{\alpha |x|}|x|^b=\gamma
\eeq
and
\beq
\label{eBLh}
\int_{\R^N} |h(x)| \, e ^{\alpha |x|}|x|^bdx<\infty,
\eeq
then, for every $z\in\R^N\setminus\{0\}$,
\[
\lim_{r \to \infty} \left(\int_{\R^N}g(x+r z)h(x) dx\right)
e^{\alpha |r z|} |r z|^b
=\gamma \int_{\R^N}h(x) e^{-\alpha\, \frac{x\cdot z}{|z|}}\, dx.
\]
\end{lemma}

\begin{lemma} \label{lem:tau-sigma}
Let $z\in\Sigma$.
For every $t\in [0,1/2]$ and $r>0$, let us set
\begin{eqnarray}
\label{101}
\delta_t(r)&\hspace{-2mm} =&\hspace{-2mm}
\left(r^{\frac{N-1}{2}}\,e^{2\sqrt{t^{\ns}\lambda_{\infty}}\,r}\right)^{-1},
\\
\label{102}
\tau_t(r)&\hspace{-2mm}=&\hspace{-2mm}\frac{2}{\rho^2}\int_{\R^N} w_{t\rho^2}(x-r z)\,w_{(1-t)\rho^2}(x-r e_1)\,dx,
\\
\label{103}
\sigma_t(r)&\hspace{-2mm}=&\hspace{-2mm}  \int_{\R^N} \left[w^{p-1}_{t\rho^2}(x-r z)\,w_{(1-t)\rho^2}(x-r e_1)+
w_{t\rho^2}(x-r z)\,w^{p-1}_{(1-t)\rho^2}(x-r e_1)\right]\!dx.
\end{eqnarray}
Then the following facts hold:

\begin{itemize}
\item[$1)$]  if $t\in\left[0,\frac 12 \right)$, then  $\displaystyle{\frac{\tau_t(r)}{\delta_t(r)}\to c_{1,t}:=
c_t\int_{\R^N} w_{(1-t)\rho^2}(y)\, e^{-\sqrt{t^s\lambda_{\infty}}\,\frac{y\cdot(e_1-z)}{2}}\,dy,
\quad \mbox{as $r\to \infty$}}$,

\item[$2)$]  if $t\in\left[0,\frac 12 \right]$, then $\displaystyle{\frac{\sigma_t(r)}{\delta_t(r)}\to c_{2,t}:=2
c_t\int_{\R^N} w^{p-1}_{(1-t)\rho^2}(y)\, e^{-\sqrt{t^s\lambda_{\infty}}\,\frac{y\cdot(e_1-z)}{2}}\,dy,
\quad \mbox{as $r\to \infty$}}$
\end{itemize}
where $c_t=c_1 t^{\ns}\left(\frac {1}{p-2}-\frac{N-1}{4}\right)$ (see \eqref{eqn:decrescita_w}).
Moreover,
\begin{itemize}
\item[$3)$]  $ c_{1,t}\cdot \left(\frac{1}{2}-t\right) \leq C$ for every $t\in\left[0,\frac 12\right)$.
 \end{itemize}
\end{lemma}

\medskip

\begin{rem}
The definition of $\tau_t(r),\sigma_t(r)$ is independent of $z$,\ by symmetry.
Moreover,  $c_{1,t}\to\infty$ as $t\to\frac 12$ and, clearly, $\tau_0(r)=\sigma_0(r)\equiv 0$.
\end{rem}

\proof
Assertions $1),\,2)$ easily follow using (\ref{eqn:decrescita_w_k}) and by Lemma \ref{BL}.
Let us prove assertion $3)$.

Without loss of generality we may assume $z=-e_1$, so that $\frac{z-e_1}{2}=-e_1$, and to simplify the notations we consider $\lambda_\infty=1$.
Moreover, since $t\mapsto c_{1,t}$ is a continuous function on $[0,1/2]$, it is sufficient to prove
$$
\lim_{t\to \frac12}\left(\frac{1}{2}-t\right)\,
\int_{\R^N} w_{(1-t)\rho^2}(y)\, e^{t^{s/2}\,y_1}\,dy
<\infty,
$$
where $y_1$ is  the first component of $y\in\R^N$.
Moreover, by (\ref{eqn:decrescita_w_k}) we are left to prove that
\beq
\label{800}
\begin{split}
&\left(\frac{1}{2}-t\right)\,
\int_{\R^N\setminus B_1(0)} |y|^{-\frac{N-1}{2}}\,e^{-(1-t)^{s/2}|y|}\, e^{t^{s/2}\,y_1}\,dy=\\
&=\left(\frac{1}{2}-t\right)\,
\int_{\R^N\setminus B_1(0)} |y|^{-\frac{N-1}{2}}\,e^{-\left[(1-t)^{s/2}-t^{s/2}\right]|y|}\, e^{-t^{s/2}(|y|-y_1)}\,dy
\le C <\infty
\quad \mbox{as $t\to (1/2)^-$}.
\end{split}
\eeq

If we set $a=\frac{1}{2}-t$, then $a\to 0^+$ as $t\to (1/2)^-$, and
since $\left[(1-t)^{s/2}-t^{s/2}\right]\ge ca$ as $a\to 0$, with $c>0$,
\eqref{800} can be estimated by
$$
a\,\int_{\R^N\setminus B_1(0)} |y|^{-\frac{N-1}{2}}\,e^{-ca|y|}\, e^{-\left(\frac{1}{2}-a\right)^{s/2}(|y|-y_1)}\,dy\le
a\,\int_{\R^N\setminus B_1(0)} |y|^{-\frac{N-1}{2}}\,e^{-ca|y|}\, e^{-\left(\frac 14\right)^{s/2}(|y|-y_1)}\,dy.
$$

Making use of spherical coordinates in the subspace
$e_1^{\perp}=\left\{v\in\R^N\ :\ v\cdot e_1=0\right\}$, denoting by $b=\left(\frac 14\right)^{s/2}$ and $\R^2_+=\R\times[0,+\infty)$,
we have to analize
$$
a\,\int_{\R^2_+\setminus B_1(0)} \frac{r^{N-2}}{\left(r^2+y_1^2\right)^{\frac{N-1}{4}}
\,e^{ca\left(r^2+y_1^2\right)^{1/2}}\, e^{b\left[\left(r^2+y_1^2\right)^{1/2}-y_1\right]}}\,dr\,dy_1=
$$
$$
=a\int_0^{\pi}\left(\int_1^{+\infty} \frac{\rho^{(N-1)/2}}{e^{[ca+b(1-\cos{\vartheta})]\rho}}\,d\rho\right)\,
(\sin{\vartheta})^{N-2}\,d\vartheta.
$$
Setting $k=ca+b(1-\cos{\vartheta})$, we can estimate
$$
\int_1^{+\infty} \frac{\rho^{(N-1)/2}}{e^{[ca+b(1-\cos{\vartheta})]\rho}}\,d\rho=
\frac{1}{k^{(N-1)/2}}\,\int_1^{+\infty} \frac{(k\rho)^{(N-1)/2}}{e^{k\rho}}\,d\rho\le
\frac{1}{k^{(N+1)/2}}\,\int_0^{+\infty} \frac{\mu^{(N-1)/2}}{e^{\mu}}\,d\mu.
$$
Then,
\begin{equation*}
\begin{split}
a\int_0^{\pi}\left(\int_1^{+\infty} \frac{\rho^{(N-1)/2}}{e^{[ca+b(1-\cos{\vartheta})]\rho}}\,d\rho\right)\,
(\sin{\vartheta})^{N-2}\,d\vartheta&\leq
C a\int_0^{\pi} \frac{(\sin{\vartheta})^{N-2}}{[ca+b(1-\cos{\vartheta})]^{(N+1)/2}}\,d\vartheta.
\end{split}
\end{equation*}
Since $1-\cos{\vartheta}\ge\frac{1}{4}\sin^2{\vartheta}$ for small $\vartheta$,
if we set $u=\sin{\vartheta}$ and $\beta=b/(4c)$, then it is sufficient to evaluate
\begin{equation*}
\begin{split}
a\int_0^{\frac 12} \frac{u^{N-2}}{[a+\beta u^2]^{(N+1)/2}}\cdot
\frac{1}{(1-u^2)^{1/2}}\,du
&\leq \widehat{C} a^{- {1}/{2}}
\int_0^{\frac 12} \frac{\left(\sqrt{\beta/a}\,u\right)^{N-2}}
{\left[1+\left(\sqrt{\beta/a}\,u\right)^2\right]^{(N+1)/2}}\,du\\
&=\widetilde{C}  \int_0^{\frac 12\,\sqrt{\frac \beta a}}
\frac{t^{N-2}}{\left(1+t^2\right)^{(N+1)/2}}\,dt\to C\quad\mbox{ as }a\to 0,
\end{split}
\end{equation*}
and the assertion follows.

\qed

\bigskip

{\noindent{\bf Proof of Proposition \ref{lem:A_varepsilon<2m}}}\hspace{2mm}
In this proof we shall consider $R>1$ fixed such that
$\R^N\setminus\Omega\subset B_{R-1}(0)$.

In order to simplify the notations, we often omit $t,z$ and
write $\psi_r=\psi_r[t,z]$, $\delta(r)=\delta_t(r)$, $\sigma(r)=\sigma_t(r)$, $\tau(r)=\tau_t(r)$ (see \eqref{100}, \eqref{101}, \eqref{102}, \eqref{103}).

We have that
$$
E(\psi_r)=\frac{1}{2}\int_{\R^N} \left[|\nabla\psi_r|^2+V(x)\psi_r^2\right] dx -\frac{1}{p}\int_{\R^N}|\psi_r|^pdx .
$$

So, to get the statement of the Lemma, we need to estimate this two integrals.

\bigskip

Let us consider $0\leq t\leq 1/2$. In an entirely analogous way we may treat also the case $1/2<t\leq 1$.

\bigskip

Let us set:
$$
w_1=w_{t\rho^2}, \quad w_2=w_{(1-t)\rho^2}, \quad \lambda_1=\lambda_{t\rho^2}=t^{\ns}\lambda_{\infty}, \quad \lambda_2=\lambda_{(1-t)\rho^2}=(1-t)^{\ns}\lambda_{\infty},
$$
and for any $i=1,2$
$$
\quad A_i=|\nabla w_i|_2^2, \quad \quad B_i=|w_i|_p^p.
$$
Recall that
$$
-\Delta w_i+\lambda_i w_i=w_i^{p-1},
$$
namely
\begin{equation}\label{equation_i}
\int_{\R^N} \nabla w_i(x)\cdot \nabla v(x)\,dx+ \lambda_i\int_{\R^N} w_i(x)\,v(x)\,dx=
\int_{\R^N} w_i^{p-1}(x)\,v(x)\,dx,
\quad \forall v\in H^1(\R^N),
\end{equation}
and so
\begin{equation} \label{w_i}
A_i-B_i=-\lambda_i|w_i|_2^2.
\end{equation}

Evidently for $t=\frac{1}{2}$ we have $w_1=w_2$, $\lambda_1=\lambda_2$, $A_1=A_2$, $B_1=B_2$.
Moreover, we recall that
\begin{equation}\label{minimi}
\frac{1}{2}A_1-\frac{1}{p}B_1=m_{t\rho^2}=t^{1+\ns}m,
\qquad
\frac{1}{2}A_2-\frac{1}{p}B_2=m_{(1-t)\rho^2}=(1-t)^{1+\ns}m.
\end{equation}

With these notations,
$$
\psi_r(x)= \rho\,\frac{\vartheta(x)\left[w_1(x-r z)+w_2(x-r e_1)\right]}
{\left|\vartheta(\cdot)\left[w_1(\cdot-r z)+w_2(\cdot-r e_1)\right]\right|_2}.
$$

\pagebreak

\no {\underline {Estimate of $\displaystyle{\left|\vartheta(\cdot)\left[w_1(\cdot-r z)+w_2(\cdot-r e_1)\right]\right|_2}$}.

\bigskip

From above
\beq
\label{eqn:norma_2_above}
\begin{split}
\left|\vartheta(\cdot)\left[w_1(\cdot-r z)+w_2(\cdot-r e_1)\right]\right|_2^2
&\leq
\left|w_1(\cdot-r z)+w_2(\cdot-r e_1)\right|_2^2
\\
&=|w_1|_2^2+|w_2|_2^2
+2\int_{\R^N} w_1(x-r z)\,w_2(x-r e_1)\,dx
\\
&=\rho^2+2\int_{\R^N} w_1(x-r z)\,w_2(x-r e_1)\,dx
\\
&=\rho^2(1+\tau(r)).
\end{split}
\end{equation}

\medskip

From below,
\begin{equation}
\label{105}
\begin{split}
\left|\vartheta(\cdot)\left[w_1(\cdot-r z)+w_2(\cdot-r e_1)\right]\right|_2^2&=
\int_{\R^N} \left|\vartheta(x)\left[w_1(x-r z)+w_2(x-r e_1)\right]\right|^2\,dx\\
&\geq \int_{\R^N} \Big|w_1(x-r z)+w_2(x-r e_1)\Big|^2\,dx\\
&\phantom{\ge}-\int_{B_{R}(0)} \Big|w_1(x-r z)+w_2(x-r e_1)\Big|^2\,dx.
\end{split}
\end{equation}

By the asymptotic behavior of $w_1$ and $w_2$ (see \eqref{eqn:decrescita_w_k}), for any $q\geq 2$ we have
\begin{equation} \label{eqn:norma_q}
\int_{B_{R}(0)} \Big|w_1(x-r z)+w_2(x-r e_1)\Big|^q\,dx
\leq 2^{q-1} \int_{B_{R}(0)} \left[|w_1(x-r z)|^q+|w_2(x-r e_1|^q\right]\,dx=o(\delta(r)).
\end{equation}

Therefore, by \eqref{eqn:norma_2_above}, \eqref{105} and \eqref{eqn:norma_q} we get
\begin{equation} \label{eqn:norma_2_below}
\left|\vartheta(\cdot)\left[w_1(\cdot-r z)+w_2(\cdot-r e_1)\right]\right|_2^2\geq
\rho^2(1+\tau(r))+o(\delta(r)).
\end{equation}

\bigskip

\no {\underline {Estimate of $\displaystyle{\int_{\R^N} \left[|\nabla\psi_r|^2+V(x)\psi_r^2\right]\,dx}$}}.

\medskip

Now, let us estimate

$$
\hspace{-7cm} \int_{\R^N} \left|\nabla\Big(\vartheta(x)\left[w_1(x-r z)+w_2(x-r e_1)\right]\Big)\right|^2\,dx
$$
\beq
 \label{1132}
\begin{split}
\le&\int_{\R^N} \left(\left|\nabla\left[w_1(x-r z)+w_2(x-r e_1)\right]\right|^2 \right.
 \\
&  +2  \int_{\R^N} \Big[\vartheta(x) \nabla\vartheta(x)\Big]\cdot
\Big(\left[w_1(x-r z)+w_2(x-r e_1)\right] 
   \nabla \left[w_1(x-r z)+w_2(x-r e_1)\right]\Big)\,dx
\\
 & + \int_{\R^N} |\nabla \vartheta(x)|^2\,\left[w_1(x-r z)+w_2(x-r e_1)\right]^2\,dx.
\end{split}
\eeq

By direct computation, (\ref{equation_i}) and \eqref{w_i}, we obtain
\beq
\label{eqn:norma_pp1}
\begin{split}
&\hspace{-2cm} \int_{\R^N} \left|\nabla\left[w_1(x-r z)+w_2(x-r e_1)\right]\right|^2\,dx \\
 =&\left|\nabla w_1\right|_2^2+\left|\nabla w_2\right|_2^2
+2\int_{\R^N} \nabla w_1(x-r z)\cdot \nabla w_2(x-r e_1)\,dx
\\
=&\left|\nabla w_1\right|_2^2+\left|\nabla w_2\right|_2^2
-(\lambda_1+\lambda_2)\int_{\R^N} w_1(x-r z)\,w_2(x-r e_1)\,dx
\\
&+\int_{\R^N} \left[w^{p-1}_1(x-r z)\,w_2(x-r e_1)+w_1(x-r z)\,w_2^{p-1}(x-r e_1)\right]\,dx
\\
=&A_1+A_2-(\lambda_1+\lambda_2)\frac{\rho^2}{2}\tau(r)+\sigma(r).
\end{split}
\eeq

\medskip

Since $\D\vartheta\equiv 0$ on $\R^N\setminus B_R(0)$, by using \eqref{eqn:decrescita_w_k} \eqref{eqn:decrescita_w_k'} we get

\begin{equation} \label{eqn:norma_pp3}
\begin{split}
\int_{\R^N} |\nabla \vartheta(x)|^2\,\left[w_1(x-r z)+w_2(x-r e_1)\right]^2\,dx
&\le
c \int_{B_R(0)}  \left[w_1(x-r z)+w_2(x-r e_1)\right]^2\,dx\\
&=o(\delta(r))
\end{split}
\end{equation}
and
\begin{equation}
\label{eqn:norma_pp4}
\begin{split}
& \left|\int_{\R^N}  \Big[\vartheta(x) \nabla\vartheta(x)\Big]\cdot
\Big(\left[w_1(x-r z)+w_2(x-r e_1)\right]  \nabla \left[w_1(x-r z)+w_2(x-r e_1)\right]\Big)\,dx\,\right|\\
\le &
c \int_{B_R(0)} \Big|\left[w_1(x-r z)+w_2(x-r e_1)\right]   \nabla \left[w_1(x-r z)+w_2(x-r e_1)\right]\Big|\, dx\\
=& o(\delta(r)).
\end{split}
\end{equation}

According to the contribution of the potential, by (\ref{HVas}), (\ref{eqn:decrescita_w_k}) and by  Lemma \ref{BL} we have
$$
\hspace{-2cm}\int_{\R^N} V(x)\left|\,\vartheta(x)\left[w_1(x-r z)+w_2(x-r e_1)\right]\right|^2\,dx
$$
\begin{equation} \label{eqn:norma_pp5}
\leq \int_{ \R^N} V(x)\Big[w_1(x-r z)+w_2(x-r e_1)\Big]^2\,dx
=o(\delta(r)).
\end{equation}

By (\ref{eqn:norma_2_below})--(\ref{eqn:norma_pp5}) we deduce

\begin{equation} \label{eqn:norma_psi}
\int_{\R^N} \left[|\nabla\psi_r|^2+V(x)\psi_r^2\right]\,dx\leq \frac{1}{1+\tau(r)+o(\delta(r))}
\left[A_1+A_2-(\lambda_1+\lambda_2)\frac{\rho^2}{2}\tau(r)+\sigma(r)+o(\delta(r))\right].
\end{equation}

\bigskip

{\underline {Estimate of $\int_{\R^N}|\psi_r|^pdx$}}.

\bigskip

Since $0\le \vartheta(x)\le 1$ in $\R^N$ and $\vartheta\equiv 1$ in
$\R^N\setminus B_{R}(0)$,  by \eqref{eqn:norma_q} and by Lemma \ref{L4.4},
we get
$$
\left|\vartheta(\cdot)\left[w_1(\cdot-r z)+w_2(\cdot-r e_1)\right]\right|_p^p =
\int_{\R^N} \Big|\vartheta(x)\left[w_1(x-r z)+w_2(x-r e_1)\right]\Big|^p\,dx  $$
\begin{equation}
\label{eqn:norma p}
\begin{split}
\geq& \int_{\R^N} \Big|w_1(x-r z)+w_2(x-r e_1)\Big|^p\,dx-\int_{B_{R}(0)} \Big|w_1(x-r z)+w_2(x-r e_1)\Big|^p\,dx \\
\geq &\left|w_1\right|_p^p+\left|w_2\right|_p^p\\
& + (p-1)\int_{\R^N} \left[w^{p-1}_1(x-r z)\,w_2(x-r e_1)+w_1(x-r z)\,w_2^{p-1}(x-r e_1)\right]\,dx\\
&-\int_{B_{R}(0)} \Big|w_1(x-r z)+w_2(x-r e_1)\Big|^p\,dx\\
=& B_1+B_2+(p-1)\sigma(r)-\int_{B_{R}(0)} \Big|w_1(x-r z)+w_2(x-r e_1)\Big|^p\,dx
\\
\geq & B_1+B_2+(p-1)\sigma(r)+o(\delta(r)).
\end{split}
\end{equation}

Taking into account (\ref{eqn:norma_2_above}) and (\ref{eqn:norma p}) we have that
\begin{equation} \label{eqn:norma lp_psi}
|\psi_r|_p^p \geq \left(\frac{1}{1+\tau(r)}\right)^{p/2}\left[B_1+B_2+(p-1)\sigma(r)+o(\delta(r))\right].
\end{equation}

\bigskip

{\underline {Estimate of $E(\psi_r)$}}.

\medskip

Therefore
\begin{equation*}
\begin{split}
E(\psi_r)=&\frac{1}{2}\int_{\R^N} \left[|\nabla\psi_r|^2+V(x)\psi_r^2\right]\,dx -\frac{1}{p}|\psi_r|_p^p\leq
\\
\leq &\frac{1}{2}\left(\frac{1}{1+\tau(r)+o(\delta(r))}\right)
\left[A_1+A_2-(\lambda_1+\lambda_2)\frac{\rho^2}{2}\tau(r)+\sigma(r)+o(\delta(r))\right]
\\
  &
-\frac{1}{p}\left(\frac{1}{1+\tau(r)}\right)^{p/2}\left[B_1+B_2+(p-1)\sigma(r)+o(\delta(r))\right].
\end{split}
\end{equation*}

\bigskip

Observe that for every $\sigma,\tau>0$ and $\delta\to 0$ we have
\begin{eqnarray*}
\frac{1}{2}\left(\frac{1}{1+\tau+o(\delta)}\right)
\left[A_1+A_2-(\lambda_1+\lambda_2)\frac{\rho^2}{2}\tau+\sigma+o(\delta)\right] && \\
 -\frac{1}{p}\left(\frac{1}{1+\tau}\right)^{\frac p 2}\left[B_1+B_2+(p-1)\sigma+o(\delta )\right]& =&\varphi(\tau,\sigma)+o(\delta ),
\end{eqnarray*}
where
$$
\varphi(\tau,\sigma)=\frac{1}{2}\left(\frac{1}{1+\tau}\right)\left[A_1+A_2-(\lambda_1+\lambda_2)\frac{\rho^2}{2}\tau+\sigma\right]
-\frac{1}{p}\left(\frac{1}{1+\tau}\right)^{\frac p2}\left[B_1+B_2+(p-1)\sigma\right],
$$
and we have performed a Taylor expansion with respect to $o(\delta )$.
We will write $\varphi_t$ in place of $\varphi$ when we want to emphasize the role of $t$.

Now, we are going to consider the Taylor  expansion of $\varphi$.
Observe that these expansions are consistent by (1) and (2) in Lemma \ref{lem:tau-sigma}.

\medskip

By (\ref{minimi}), and taking into account  (\ref{w_i}), we have that
\begin{eqnarray}
\label{300}
\varphi(0,0)&=&\left(\frac{1}{2}A_1-\frac{1}{p}B_1\right)+\left(\frac{1}{2}A_2-\frac{1}{p}B_2\right)
=\left[t^{1+\ns}+(1-t)^{1+\ns}\right]m,
\\ \noalign{\bigskip}
\nonumber \displaystyle{\frac{\partial \varphi}{\partial \tau}(0,0)} &=&
\displaystyle{-\frac{1}{2}\left[A_1+A_2+(\lambda_1+\lambda_2)\frac{\rho^2}{2}-B_1-B_2\right]} \\
\noalign{\medskip} \nonumber
&= &\displaystyle{-\frac{1}{2}\left[A_1-B_1+\lambda_1\frac{\rho^2}{2}\right]
-\frac{1}{2}\left[A_2-B_2+\lambda_2\frac{\rho^2}{2}\right]} \\
\noalign{\medskip}
\label{301}&= &\displaystyle{-\frac{1}{2}\left[-\lambda_1 t\rho^2+\lambda_1\frac{\rho^2}{2}\right]
-\frac{1}{2}\left[-\lambda_2 (1-t)\rho^2+\lambda_2\frac{\rho^2}{2}\right]} \\
\noalign{\medskip}
\nonumber &= &\displaystyle{-\frac{1}{2}\lambda_1 \rho^2\left(\frac{1}{2}-t\right)
+\frac{1}{2}\lambda_2 \rho^2\left(\frac{1}{2}-t\right)} \\
\noalign{\medskip}
\nonumber &= &\displaystyle{\frac{1}{2}\left(\frac{1}{2}-t\right)\left(\lambda_2-\lambda_1\right) \rho^2}. \\
 \noalign{\bigskip}
\label{302} \frac{\partial \varphi}{\partial \sigma}(0,0) &= &
-\left(\frac{1}{2}-\frac{1}{p}\right).
\end{eqnarray}

Hence, we obtain
$$
\varphi(\tau,\sigma)= \varphi(0,0)+\frac{1}{2}\left(\frac{1}{2}-t\right)\left(\lambda_2-\lambda_1\right) \rho^2\tau
-\left(\frac{1}{2}-\frac{1}{p}\right)\sigma+o\left(\sqrt{\tau^2+\sigma^2}\right),
\quad \mbox{as $(\tau,\sigma)\to (0,0)$}.
$$

Now, we want to analyze the asymptotic behaviour of  $\varphi(\tau(r),\sigma(r))$ and, as a consequence, of $E(\psi_r)$, as $r\to\infty$.

By Lemma \ref{lem:tau-sigma},  for every $t\in\left[0,\frac 12\right)$   we have
$$
\varphi(\tau_t(r),\sigma_t(r))= \varphi(0,0)+\left[\frac{1}{2}\left(\frac{1}{2}-t\right)\left(\lambda_2-\lambda_1\right) \rho^2c_{1,t}
-\left(\frac{1}{2}-\frac{1}{p}\right)c_{2,t}\right]\delta_t(r)+o(\delta_t(r)),
$$
as  $r\to \infty$.
Taking into account Lemma \ref{lem:tau-sigma} and $\lambda_1\to\lambda_2$ as $t\to\frac 12$, we get
$$
\left[
\frac{1}{2}\left(\frac{1}{2}-t\right)\left(\lambda_2-\lambda_1\right) \rho^2c_{1,t}
-\left(\frac{1}{2}-\frac{1}{p}\right)c_{2,t}\right]\lo
-\left(\frac 12 -\frac 1p\right) C_2<0,\qquad\mbox{ as }t\to\left(\frac 12\right)^-.
$$

Then, for suitable constants $C>0$ and $\mu\in \left( 0,\frac 12\right)$
$$
E(\psi_r[t,z])\le [t^{1+s}+(1-t)^{1+s}]\, m-C\delta_t(r)+o(\delta_t(r))\quad\forall t\in\left[\mu,  1/2\right)
$$
(see \eqref{300}}.
Since $\delta_t(r)\to 0$ as $r\to \infty$, uniformly in $t$ (see \eqref{101}), we can conclude that there exists $r_1>0$ such that
\beq
\label{900}
E(\psi_r[t,z])<\max_{t\in [0,1]}[t^{1+s}+(1-t)^{1+s}]\, m=2^{-s}m\qquad\forall t\in\left[\mu, 1/2\right),\quad \forall r>r_1.
\eeq
If $t=\frac 12$,   for suitable $\alpha,\beta,\gamma\in\R$ we infer
\beq
\label{901}
\varphi(\tau,\sigma)=
\varphi(0,0)-\left(\frac{1}{2}-\frac{1}{p}\right)\sigma+
\alpha\tau^2+\beta\sigma^2+\gamma\sigma\tau+o\left(\tau^2+\sigma^2\right)
\quad \mbox{as $(\tau,\sigma)\to (0,0)$}.
\eeq
Now, consider that by Lemma \ref{lem:tau-sigma}
\beq
\label{902}
\sigma_{1/2}(r)=(c_{2,\frac 12}+o(1))\delta_{1/2}(r)
\eeq
and that, fixed  $\eta\in\left(0,\frac{\sqrt{2^{-\ns}\lambda_{\infty}}}{2}\right]$, by Lemma \ref{BL}
$$
\tau_{1/2}(r)=o\left(e^{-2\left(\sqrt{2^{-\ns}\lambda_{\infty}}\,-\eta\right) r}\,r^{-\frac{N-1}{2}}\right),
\qquad\mbox{ as } r\to +\infty,
$$
so that
\beq
\label{903}
\left(\tau_{1/2}(r)\,\right)^2=o\left(\delta_{1/2}(r)\right).
\eeq
By \eqref{901}, \eqref{902}, \eqref{903}, we can conclude that there exists $r_2>0$ such that
\beq
\label{433}
E(\psi_r[1/2,z])=\varphi(0,0)-C\left(\frac 12-\frac 1p\right)\delta_{1/2}(r)+o\left(\delta_{1/2}(r)\right)<2^{-s}m\qquad\forall r>r_2.
\eeq

If $t\in[0,\mu]$, then, taking into account that $\delta_t(r)\le\delta_0(r)$ for every $t\in[0,1/2]$, we get
\begin{equation*}
\begin{split}
E(\psi_r[t,z]) &\le \varphi_t(0,0)+O\left(\delta_0(r)\right)\\
&\le [t^{1+s}+(1-t)^{1+s}]m+O\left(\delta_0(r)\right)\\
&\le [\mu^{1+s}+(1-\mu)^{1+s}]m+O\left(\delta_0(r)\right)\qquad \forall t \in[0,\mu].
\end{split}
\end{equation*}
Since $[\mu^{1+s}+(1-\mu)^{1+s}]m<2^{-s}m$, there exists $r_3>0$ such that
\beq
\label{904}
E(\psi_r[t,z])<2^{-s}m\qquad\forall t\in[0,\mu],\quad\forall r>r_3.
\eeq
So, assertion $(a)$ follows from \eqref{900}, \eqref{433}, \eqref{904}, for every $r>\bar r:=\max\{r_1,r_2,r_3\}$.

\bigskip

The estimates developed above also show that
$$
E(\psi_r[0,z])\le\varphi_0(0,0)+o\left(\delta_0(r)\right)=m+o\left(\delta_0(r)\right)
$$
so that $ E(\psi_r[0,z])\lo m$, as $r\to\infty$.
The same arguments work to evaluate
$$
E(\psi_r[1,z])\le m+o\left(\delta_0(r)\right),
$$
uniformly in $z\in\Sigma$.
So, also $(b)$ is proved.

\qed

\subsection{Other estimates}

In this subsection we consider the nonautonomous case $\Omega\neq\R^N$ or $V\not\equiv 0$.

\smallskip

The following   {definition of barycenter} of a function $u \in
H^1(\R^N)\setminus \left\lbrace 0 \right\rbrace, $ has been
introduced in
\cite{CP03}.
We set
\beq
\label{mu}
\mu(u)(x)=\frac{1}{|B_1(0)|}\int_{B_1(x)}|u(y)|dy
\eeq
and we remark that  $\mu(u)$ is bounded and continuous, so  we can
introduce the function
\beq
\label{hat}
\hat{u}(x)=\left[\mu(u)(x)-\frac{1}{2}\max \mu(u)\right]^{+},
\eeq
that is continuous and has compact support.
Thus, we can set
$\beta: H^1(\R^N)\setminus\{0\}\rightarrow \mathbb R^N$
as
$$
\beta(u)=\frac{1}{|\hat{u}|_{1}}\int_{\mathbb R^N}\hat{u}(x)\, x\, dx.
$$
The map $\beta$ has  the following properties:
 \begin{eqnarray}
& & \beta \mbox{ is continuous in }H^1(\mathbb R^N)\setminus \{0\};
\label{b1}
\\
&& \mbox{if } u\mbox{ is a radial function, then }\beta(u)=0;
\label{b2}
\\
&& \beta(tu)=\beta(u)
\qquad \forall t\in\R\setminus\{0\},\quad \forall u\in H^1(\mathbb R^N)\setminus
\{0\};
\label{b3}
\\
&& \beta(u(x-z))=\beta(u)+z
\qquad \forall z\in\R^N\quad \forall u\in H^1(\mathbb R^N)\setminus
\{0\}.
\label{b4}
\end{eqnarray}

Let us set
$$
C_0=\inf\{E(u):\enskip u\in H^1_0(\Omega),\,\, |u|_2=\rho,\,\, \beta(u)=0\}.
$$

\begin{lemma} \label{lem:c_0}
We have that $C_0>m$.
\end{lemma}

\proof
Of course we have that $C_0\geq m$.
Assume by contradiction that $C_0=m$.
Then by Corollary \ref{Cmin}
there exists a sequence $(y_n)_n$ in $\R^N$ such that $|y_n|\lo \infty$, as $n\to\infty$, and
$$
u_n(x)=w(x-y_n)+\phi_n(x),\qquad \phi_n\to 0\ \mbox{ strongly in }H^1(\R^N).
$$
By (\ref{b1}), \eqref{b2} and  (\ref{b4}) we have
$$
0=\beta(u_n)=\beta(w(\cdot-y_n)+\phi_n)=
\beta(w+\phi(\cdot+y_n))+y_n=y_n+o(1),
$$
contrary to $|y_n|\lo \infty$, as $n\to\infty$.

So the proof is completed

\qed

\begin{lemma} \label{lem:c_0_A_varepsilon}
Let ${\cal A}_r$ and $\cL_r$ be as in Proposition \ref{lem:A_varepsilon<2m}.
Then $\widehat{r}>0$ exists such that
\beq
\label{Me}
 \cL_r<C_0\leq {\cal A}_r,\qquad\mbox{ for all }r\geq \widehat{r}.
\eeq
\end{lemma}

\proof
 Inequality $\cL_r<C_0$, for large $r$, follows from \eqref{936}  and Lemma \ref{lem:c_0}.

\smallskip

To get the second inequality in \eqref{Me},
we claim that
$\beta\big(\vartheta(\cdot)w(\cdot-rz)\big)\cdot z >0$ for all $z\in\Sigma$, for large $r$.
Indeed, by (\ref{b1}),\eqref{b2} and (\ref{b4}) we have
$$
\Big|\beta\big(\vartheta(\cdot)w(\cdot-rz)\big)-rz\Big|=
\Big|\beta\big(\vartheta(\cdot+rz)\, w\big)\Big|\xrightarrow[]{r\rightarrow \infty}0,
$$
because $\vartheta(\cdot+rz)w\to w$ in $H^1$ as $r\to \infty$, by \eqref{eqn:decrescita_w} and \eqref{eqn:decrescita_w'}.
Hence
$$
\beta(\vartheta(\cdot)w(\cdot-rz))=r z+o(1)\qquad\mbox{ as }r\to\infty,
$$
as asserted.
So, for $r$ large, the deformation $\cG:[0,1]\times\Sigma\to
\R^N\setminus\{0\}$ given by
\beq
\label{1552}
\cG(s,z)=s\beta(\psi_r[1,z])+(1-s)\,z
\eeq
is well defined.
Then we claim that there exists $(t_r,z_r)\in [0,1]\times\Sigma$ such that
\beq
\label{td}
\beta(\psi_r[t_r,z_r])=0.
\eeq
Indeed, by the continuity of the maps $\beta$ and $\psi_r$, by  $\cG(s,z)\neq 0$, $\forall (s,z)\in[0,1]\times\Sigma$, and by the invariance of the topological degree by homotopy, we have
$$
0\neq d(\Id,\Sigma\times[0,1),0)=d(\beta\circ\psi_r,\Sigma\times[0,1),0).
$$
From \eqref{td} there follows $C_0\leq E(\psi_r[t_r,z_r])\leq {\cal A}_r$, that completes the proof.

\qed


\sezione{Proof of the main results}

The existence of a positive solution for the autonomous case $\Omega=\R^N$ and $V\equiv 0$ is well known, so we prove our results when  $\Omega=\R^N$ and $V\equiv 0$ does not occur.

\medskip

{\noindent{\bf Proof of Theorem \ref{T}}}\hspace{2mm}
Let us recall the values
\beq
\label{1201b}
{\cal A}_r=\max\left\{E\left(\psi_{r}[t,z]\right):\enskip t\in [0,1],\,\, z\in\Sigma\right\},
\eeq
\beq
\label{1202b}
{\cal L}_r=\max\{E\left(\psi_r[1,y]\right):\enskip z\in \Sigma\},
\eeq
\beq
\label{1644}
C_0=\inf\{E(u):\enskip u\in H^1_0(\Omega),\,\, |u|_2=\rho,\,\, \beta(u)=0\}.
\eeq
By Propositions \ref{lem:A_varepsilon<2m} and
\ref{prop:no_ground_state} (see \eqref{1229}), and by Lemma  \ref{lem:c_0_A_varepsilon}, we have that
for all $r>\max\{\bar{r},\widehat{r}\}$
\beq
\label{1635}
m<{\cal L}_r<C_0\leq {\cal A}_r <2^{-s}m.
\eeq

We {\em  claim} that the functional $E$ has a (PS)-sequence in
$[C_0,{\cal A}_r]$.
This done, the existence of a critical point $\bar u$ of $E$ on $S_{\rho}$ with
$E(\bar u)\le \cA_r$ follows from Proposition \ref{PPS}.

Assume, by contradiction, that no (PS)-sequence exists in
$[C_0,{\cal A}_r]$.
Then,  usual deformation arguments imply the existence of $\eta>0$ such that
the sublevel
$E^{C_0-\eta}:=\{u\in H^1_0(\Omega):\enskip |u|_2^2=\rho^2,\,\, E(u)\leq
C_0-\eta\}$ is a deformation retract of the sublevel
$E^{{\cal A}_r}:=\{u\in H^1_0(\Omega):\enskip |u|_2^2=\rho^2,\,\, E(u)\leq {\cal A}_r\}$,
namely there exists a continuous function
$\varphi:E^{{\cal A}_r}\to E^{C_0-\eta}$ such that
\beq
\label{1553}
\varphi(u)=u\qquad\mbox{ for any }u\in E^{C_0-\eta}.
\eeq
Furthermore, by (\ref{1635}) we can also assume $\eta$ so small that
\beq
\label{1603}
C_0-\eta>{\cal L}_r.
\eeq
Let us define the map $\cH:[0,1]\times\Sigma\to\R^N$ by
$$
\cH(s,y)=\beta\left(\varphi\big(\psi_r[t,z]\big)\right).
$$
By (\ref{1603}),  (\ref{1553}), and by using the map $\cG$ introduced
in (\ref{1552}), we deduce that $\cH$ maps
$\{1\}\times  \Sigma$ in a set homotopically equivalent to
$ \Sigma$  in $\R^N\setminus\{0\}$.
Since $\cH$ is a continuous map, and arguing exactly as for \eqref{td}, we get the existence of a point $(\tilde t,\tilde z)\in [0,1]\times \Sigma$ such that
$$
0=\cH(\tilde t,\tilde z)=\beta(\varphi(\psi_r[\tilde t,\tilde z])).
$$
Then by the definition of $C_0$ we see $E(\varphi(\psi_r[\tilde t,\tilde z]))\ge C_0$, contrary to
$\varphi\big(\psi_r[t,z]\big)\in
E^{C_0-\eta}$ for every $(t,z)\in[0,1]\times \Sigma$, so the
claim must be true.

\smallskip

Finally, since $E(\bar u)\in (m,2^{-s}m)$ then $\bar u$ has constant sign by Proposition \ref{Psegno}.
Observe that since $\bar u$ solves \eqref{P} if and only if $-\bar u$ solves \eqref{P}, then we have a nonnegative solution.

\qed

\medskip

{\noindent{\bf Proof of Theorem \ref{T2}}}\hspace{2mm}
This proof proceed in two different ways, according to $\Omega\neq\R^N$ or $\Omega=\R^N$.
In both cases, we identify a topological configuration analogous to the one employed in the proof of Theorem \ref{T} to get the solution.
We only outline the procedure, because the argument is the same already developed in details.

\bigskip

{\underline {Case $\Omega\neq\R^N$ }}.

\medskip

Let us apply Proposition \ref{lem:A_varepsilon<2m}, Lemma \ref{lem:c_0} and Lemma \ref{lem:c_0_A_varepsilon} with $V\equiv 0$.
Then for a fixed $r>0$, large enough,  we get
\beq
\label{1200}
m<{\cal L}_{r,0}<C_{0,0}\leq {\cal A}_{r,0} <2^{-s}m,
\eeq
where $ {\cal L}_{r,0},\, C_{0,0},\  {\cal A}_{r,0}$ are defined as in \eqref{1201b}, \eqref{1202b}, \eqref{1644}, with the functional $E$ replaced by
$$
E_0(u):=\frac 12 \int_\Omega|\D u|^2dx-\frac 1p \int_\Omega|u|^pdx,\qquad  u\in H^1_0(\Omega).
$$
The configuration in \eqref{1200} depends on $\Omega$ and $\rho$.

Now, observe that $ {\cal L}_{r,0}\le {\cal L}_{r}$, $  C_{0,0}\le C_{0}$,   ${\cal A}_{r,0}\le {\cal A}_{r}$, by \eqref{HV}. 
Moreover, by the H\"older inequality and \eqref{eqn:decrescita_w_k}, it follows that
$$
\lim_{|V|_q\to 0}\int_{\R^N}V(x)\left(\psi_{r}[t,z]\right)^2  dx\le \lim_{|V|_q\to 0}|V|_q\max_{t\in[0,1],\ z\in\Sigma}|(\psi_r[t,z])^2|_{q'}=0,
$$
uniformly in $t\in[0,1]$ and  $z\in\Sigma$, that implies
$$
\lim_{|V|_q\to0}{\cal L}_{r}= {\cal L}_{r,0},\quad   \lim_{|V|_q\to0}{\cal A}_{r}= {\cal A}_{r,0}.
$$

Hence, taking also into account Lemma  \ref{lem:c_0_A_varepsilon}, we see that there exists $L=L(\Omega,\rho)$ such that if $|V|_q<L$ then the configuration \eqref{1635} is restored.
In particular, 
\beq
\label{10}
L<\frac{1}{\max\limits_{t\in[0,1],\ z\in\Sigma}|(\psi_r[t,z])^2|_{q'}}\, (2^{-s}m-{\cal A}_{r,0}).
\eeq 
As a consequence, if  $|V|_q<L$  we get a critical value for $E$ constrained on $S_\rho$, in the energy range $(m,2^{-s}m)$.
 
 \smallskip

Again, since the solution $\bar u$ we have found verifies $E(\bar u)\in (m,2^{-s}m)$,   it is a constant sign solution by Proposition \ref{Psegno}.

\bigskip

{\underline {Case $\Omega=\R^N$ }}.

\medskip

Let $r>0$ and let us introduce the values
\beq
\label{1201}
\widetilde {\cal A}_r=\max\left\{E\left(w(\cdot -y)\right):\enskip y\in B_r(0)\right\} ,
\eeq
\beq
\label{1202}
\widetilde {\cal L}_r=\max\{E\left(w(\cdot -y)\right):\enskip y\in \partial B_r(0)\}.
\eeq
Then it turns out that for every fixed $V$ there exists  $r_V>0$ such that for every $r> r_V$
\beq
\label{5}
m<\widetilde {\cal L}_r<C_0\le \widetilde {\cal A}_r.
\eeq
By H\"older inequality, for every $y\in\R^N$
$$
E\big(w(\cdot -y)\big)=m+\frac 12 \int_{\R^N} V(x)\, [w(\cdot -y)]^2dx\le m+\frac 12 |V|_q\, |w^2|_{q'}.
$$
Hence, $\widetilde\cA_r<2^{-s}m$ provided
\beq
\label{11}
\frac 12 |V|_q\, |w^2|_{q'}<\left(1-\frac{1}{2^s}\right)(-m).
\eeq
After some computation, by  \eqref{700} with 1 in place of $\rho$ and $\rho^2$ in place of $k$, we see that \eqref{11} is equivalent to 
 \beq
 \label{12}
 |V|_q\, \rho^{-\frac {2s}{q}\left( q-\frac N2\right)} <c ,
 \eeq
 for a suitable constant $c>0$ depending on $N,p$ and $q$.
 Hence,  setting for example 
 \beq
 \label{Sergio}
 L=\frac12 c \, \rho^{\frac {2s}{q}\left( q-\frac N2\right)} ,
 \eeq
  if $|V|_q<L$ then the inequalities \eqref{5} can be completed as 
\beq
\label{1204}
m<\widetilde {\cal L}_r<C_0\le \widetilde {\cal A}_r<2^{-s}m,\qquad \forall r>r_V.
\eeq
Moreover, by \eqref{12} we obtain \eqref{eaas0} and \eqref{eaas}.
Finally, by \eqref{1204} we can argue as in the proof of Theorem \ref{T} to get the solution we are looking for.

\qed

\begin{rem}
\label{ras}
Let $\rho>0$ be fixed.
Arguing as in \cite[Theorem 3.1]{98DCDS}, it is possible to verify that 
$$
\lim_{r(\Omega)\to\infty}C_{0,0}= 2^{-s}m
$$
(see\eqref{7}). 
Then by $C_{0,0}\le{\cal A}_{r,0}<2^{-s}m$ and \eqref{10} we obtain
$$
\lim_{r(\Omega)\to\infty} L=0.
$$

\end{rem}


\sezione{The case N=1}


In this section we consider the 1-dimensional case 
\beq
\label{P_1}
\tag{$P_1$}
\left\{
\begin{array}{ll}
-u''+\lambda u+V(x)\, u=|u|^{p-2}u & \mbox{in $I$}, \\
 \noalign{\medskip}
\lambda\in\R, \quad u\in H^1_0(I) , \quad \displaystyle{\int_{I} u^2\,dx=\rho^2},
\end{array}
\right.
\eeq
 where we can consider $I=\R$ or $I=(0,\infty)$, $V\in L^q(I)$ for some $q\in [1,\infty]$, $V\ge 0$  and $2<p<6$.

\medskip

First, let us consider the entire case.
For $N=1$, in the asymptotic behaviour of the limit function $w$ (see \eqref{eqn:decrescita_w} and \eqref{eqn:decrescita_w'}) we cannot take advantage of the polinomial contribution provided by $|x|^{\frac {N-1}{2}}$.
As a consequence, the key estimate \eqref{A_r} in Proposition \ref{lem:A_varepsilon<2m} does not hold and Theorem \ref{T} fails.

\smallskip

On the other hand, Theorem \ref{T2} does not need Proposition \ref{lem:A_varepsilon<2m} and it can be again stated:
\begin{teo}
\label{T3}
Let   $I=\R$, $\rho>0$, $V\in L^q(\R)$, for $q\in[1,+\infty]$, and $V\ge 0$ in $\R$, then there exists a constant $L=L(\rho)>0$ such that if $|V|_q<L$ then  problem $(P_1)$  has a positive solution.
\end{teo}

{\it Sketch  of the proof.}
This theorem can be proved exactly as Theorem \ref{T2}-case $\Omega=\R^N$:
 we introduce  $\widetilde {\cal A}_r$, $C_0$ and $\widetilde {\cal L}_r$  (see \eqref{1201}, \eqref{1644} and \eqref{1202}) and then  observe that for $|V|_q$ small $m <\widetilde {\cal L}_r<C_0\le  \widetilde {\cal A}_r<2^{-s}m$.
 So, by Proposition \ref{PPS}, we can argue by deformation as in the proof of Theorem \ref{T}, getting the existence of a nonnegative solution $\bar u$.
By Proposition  \ref{stab}, the solution $\bar u$ is positive.

 \qed

According to the exterior case, that is $I=(0,\infty)$, first let us state some nonexistence results.

\begin{rem}
\label{R:EL}
If  $I=(0,+\infty)$ and $V(x)\equiv 0$, then the autonomous problem $(P_1)$ has no solution.
\end{rem}
Indeed,  if $u$ is a solution of  $(P_1)$, then it is a regular free solution of an equation of the form $-u''=f(u)$, $u(0)=0$, where $f(u)=-\lambda u+ |u|^{p-2}u$.
Then $u\equiv 0$ by \cite[Remark I.3]{EL82PRSE}, contrary to $u\in S_\rho$.

\begin{prop}
\label{P:NEN=1}
Let $I=\R$ or $I=(0,\infty)$. 
If $V\in L^\infty(I)$ is a monotone locally Lipshitz function such that $V'\neq 0$ on a set of positive measure, then problem  $(P_1)$ has no  solution.
\end{prop}

Proposition \ref{P:NEN=1} is a simplified version of Proposition \ref{P:NE} in the 1-dimensional case, adapted also to half lines.

\proof
We are assuming $V$ non increasing  and $I=(0,\infty)$.
 
Suppose by contradiction that  there exists a solution  $u$ of $(P_1)$.
Since the solutions of $(P_1)$ are in $H^2(I)$, we can consider a sequence $(u_n)_n$ in $\cC^\infty_0\big((1/n,\infty)\big)$ such that  $u_n\to u$ in $ H^2(I)$.
For every $n\in\N$, the  map $t\mapsto u_n(x+t)$, $t\in(-1/n,\infty)$, turns out to be a smooth  curve in $H^1(I)$ and we can consider the regular map $t\mapsto f_n(t):=E\big(u_n(x+t)\big)$.
Since $u$ is a solution of $(P_1)$, we have
\beq
\label{1100}
f_n'(0)=E'(u_n)[u_n']=O(1).
\eeq
From the monotonicity of $V$ and Fatou's Lemma  we infer
\beq
\label{1101}
\begin{split}
f_n'(0) & =\lim_{t\to 0}\frac 1t\, \left[ E\big(u_n(\cdot+t)\big)-E\big(u_n\big)\right]\\
& =\lim_{t\to 0^+}\frac {1}{2t}\int_{\frac 1n}^\infty V(x)[u_n^2(x+t)-u_n^2(x)]\, dx\\
& =\lim_{t\to 0^+}\frac 12\int_{\frac 1n}^\infty\frac{V(x-t)-V(x)}{t}\, u_n^2(x)\, dx\\
& \ge \frac 12\int_0^\infty (-V'(x))\, u_n^2(x)\, dx \\
&\ge 0.
\end{split}
\eeq
Letting $n\to\infty$ in \eqref{1100} and \eqref{1101}, we obtain
$$
\int_0^\infty V'(x)u^2(x)\, dx =0,
$$
that is a contradiction.
Indeed, by assumption $V'\le 0$ a.e. in $I$ and $V'<0$ on a set of positive measure, while $|\{x\in I\, :\, u(x)=0\}|=0$ because $u(x)=0$ implies $u'(x)\neq 0$ otherwise $u\equiv 0$ by the Cauchy-Lipschitz theorem.

\smallskip

If $I=\R$ we can proceed in a similar way. 
 
 \qed

\medskip

\begin{rem}
In the proof of Theorem \ref{P:NEN=1}, we cannot consider directly the curve in $H^1_0(I)$ defined by  $\gamma(t):=u(\cdot+t)$, $t\ge 0$.
Indeed, $u'(0)\neq 0$ prevents $u'\in H^1_0(I)$, so $\gamma$  would be not a regular curve in $H^1_0(I)$.
\end{rem}

\medskip

Finally, let us state an existence result on half lines, that inherits the topological structure of the entire case.

\begin{teo}
Let $I=(0,\infty)$, $\rho>0$, $V\in L^q(\R)$, for some $q\in [1,+\infty]$ and  $V\ge 0$.  
If $|V|_q<L$, with $L$ as in Theorem \ref{T3}, then there exists $\bar R=\bar R(V,\rho)>0$ such that problem $(P_1)$ with $V(x-R)$ has a positive solution for every $R>\bar R$.
\end{teo}

\proof
In order to get a solution, we are going to solve $(P_1)$ with $V$ fixed, on $H^1_0\big((-R,\infty)\big)$.

Let us consider $(P_1)$  on $H^1 (\R)$,  with $V$ fixed, and let us define $\widetilde {\cal A}_r$, $C_0$ and $\widetilde {\cal L}_r$ as in the proof of Theorem \ref{T3}.
Moreover, let us fix a cut-off function $\widetilde\vartheta\in C^\infty(\R,[0,1])$ such that $\widetilde\vartheta(x)=0$ $\forall x\in(-\infty,0)$, $\widetilde\vartheta(x)=1$ $\forall x\in(1,\infty)$ and introduce
 $$
\widetilde {\cal A}_{r,R}=\max\left\{E\left(\widetilde\vartheta(\cdot+R)w(\cdot -y)\right):\enskip y\in [ -r,r]\right\} ,
$$
$$
C_{0,R}=\inf\{E(u)\ :\ u\in H^1(-R,\infty),\ |u|_2=\rho,\  \beta (u)=0\},
$$
$$
\widetilde {\cal L}_{r,R}=\max\big\{E\left(\widetilde\vartheta(\cdot+R)w(\cdot -y)\right):\enskip y\in\{-r,r\}\big\}.
$$

Then,
$$
\lim_{R\to\infty}\widetilde {\cal A}_{r,R}=\widetilde {\cal A}_{r},\qquad
\lim_{R\to\infty}C_{0,R}=C_{0},\qquad
\lim_{R\to\infty}\widetilde {\cal L}_{r,R}=\widetilde {\cal L}_{r}.
$$

Hence, for large $R$, we get
$$
m <\widetilde {\cal L}_{r,R}<C_{0,R}\le  \widetilde {\cal A}_{r,R}<2^{-s}m
$$
and we can argue as in the proof of Theorem \ref{T}, by Proposition \ref{PPS}.

\qed


\appendix
\section*{Appendix}


\renewcommand{\theequation}{A.\arabic{equation}}
\setcounter{equation}{0}

The following proof could be shortened by taking into account, for example, Lemma B.3 in \cite{St08book}. 
We develop some more details on the bootstrap procedure, for the sake of completeness.

For classical regularity results, we refer the reader for example to  \cite[\S 10]{LiebLossBook}, and in particular Theorem 10.2 there, or to \cite[\S 8]{GiTr01book}.

\medskip

{\em Proof of Proposition \ref{stab}}.\ 
 $(a)$ \
Once proved   $\bar u\in\cC^{0,\alpha}_{\loc}(\Omega)$  we are done. 
Indeed, when $q>\frac N2$ we can apply the Harnack inequality and conclude that the nonnegative solutions are actually strictly positive in $\Omega$  (see \cite[Theorem 7.2.1]{PucciSerrinBook}). 

\smallskip

To verify that   $\bar u\in\cC^{0,\alpha}_{\loc}(\Omega)$, we begin by observing that the function $\bar u$ verifies 
\beq
\label{eqn:stab}
-\Delta \bar u=\psi(x)\bar u\quad\mbox{ for }\ \psi(x)=-\lambda  -V(x) +|\bar u|^{p-2},\qquad x\in\Omega.
\eeq

Since $V\in L^q_{\loc}(\Omega)$ with $q>N/2$, we have that
$\psi\in L^{q_1}_{\loc}(\Omega)$ with
$$
q_1=\min\left\{\frac{2N}{(N-2)(p-2)},\, q\right\}.
$$
Notice also that  $q_1>\frac{N}{2}$ because $p<2^*$.
By H\"older inequality, $\psi \bar u\in L^{r_1}_{\loc}(\Omega)$ with $r_1$ defined by
$$
\frac{1}{r_1}=\frac{1}{q_1}+\frac{N-2}{2N}.
$$
By regularity results, we have that:

\begin{itemize}

\item[$1)$] if $r_1>\frac{N}{2}$, then $\bar u\in\cC^{0,\alpha}_{\loc}(\Omega)$ with $\alpha<2-\frac{N}{r_1}$;

\item[$2)$] if $r_1=\frac N2$, then $\bar u\in L^s_{\loc}(\Omega)$ $\forall s>1$; 

\item[$3)$] if $r_1<\frac{N}{2}$, then $\bar u\in L^{r_1N/(N-2r_1)}_{\loc}(\Omega)$.

\end{itemize}

In case (1) the assertion is proved.

\smallskip

In case (2) it is readily seen that $\psi \, \bar u\in L^{s_1}_{\loc}(\Omega)$ for every $1\le s_1<q$.
Hence we can conclude as in the previous case, choosing $N/2<s_1<q$.

\smallskip

If (3) holds, then $\psi \in L^{q_2}_{\loc}(\Omega)$, with
$$
q_2=\min\left\{\frac{r_1 N}{(N-2r_1)(p-2)},\, q\right\}.
$$
Since $q_1>\frac{N}{2}$, we get $\frac{r_1 N}{N-2r_1}>\frac{2N}{N-2}$ and hence $q_2\geq q_1$.

Now, repeating the same argument of the previous step,
by H\"older inequality $\psi \bar u\in L^{r_2}_{\loc}(\R^N)$ with $r_2$ defined by
$$
\frac{1}{r_2}=\frac{1}{q_2}+\frac{N-2r_1}{r_1 N},
$$
and again we have

\begin{itemize}

\item[$4)$] if $r_2\ge \frac{N}{2}$, then $\bar u\in\cC^{0,\alpha}_{\loc}(\Omega)$;

\item[$5)$] if $r_2<\frac{N}{2}$, then 
$\bar u\in L^{r_2N/(N-2r_2)}_{\loc}(\Omega)$.

\end{itemize}

Observe that $r_2>r_1$ because $q_2\geq q_1$ and $\frac{r_1 N}{N-2r_1}>\frac{2N}{N-2}$.

Iterating this bootstrap argument, we claim that, after $\bar k$ steps, $r_{\bar k}\ge \frac N2$, so that we are done.
If, by contradiction, the claim is false, then for every integer $k$ we define 
\beq
\label{passo_k}
q_k=\min\left\{\frac{r_{k-1} N}{(N-2r_{k-1})(p-2)},\, q\right\},
\qquad
\frac{1}{r_k}=\frac{1}{q_k}+\frac{N-2r_{k-1}}{r_{k-1} N}.
\eeq
Inductively, it turns out that $q_k\geq q_{k-1}$ and $r_k>r_{k-1}$, for any $k\in \N$,
with $q_k\le q$ and $r_k<\frac N2$.
Let us set
$$
\lim_{k\to +\infty} r_k=R
\qquad
\lim_{k\to +\infty} q_k=Q.
$$
Getting $k\to  \infty$ in (\ref{passo_k}) we obtain
$$
\frac{1}{R}=\frac{1}{Q}+\frac{N-2R}{R N}
$$
that implies $Q=\frac{N}{2}$, contrary to  $Q\geq q_1>\frac{N}{2}$.
So, a contradiction arises and $(a)$ is proved.
 
\medskip
 
 {$(b)$}\ The function $\bar u$ verifies
 $$
 -\Delta \bar u=\phi \quad\mbox{ for }\ \phi(x)=-\lambda\bar u   -V(x) \bar u+|\bar u|^{p-2}\bar u,\qquad x\in\Omega.
 $$
If $V\in L^q_{\loc}(\Omega)$, with $q>N$, then $\bar u\in L^\infty_{\loc}(\Omega)$, by $(a)$, and so $\phi\in L^q_{\loc}(\Omega)$, 
that allows us to deduce $\bar u\in \cC^{1,\alpha}(\Omega)$.
 
\medskip
 
  {$(c)$} This point follows by classical regularity results.
  
\medskip
 
  {$(d)$} By classical regularity results, $\bar u$ is continuously differentiable.
  If $\bar u$ is nonnegative, then it turns out to be positive, by the Harnack inequality  
  (see \cite[Theorem 7.2.1]{PucciSerrinBook}).
    
\medskip
 
  {$(e)$} By a direct verification on  \eqref{P_1}, we see that there exists $\bar u''$ in  $L^{q}(\Omega)$, in a weak sense.
  Hence, if $q<2$ then $\bar u'\in W^{1,q}_{\loc}(\Omega)$ and if $q\ge 2$ then  $\bar u'\in H^1_{\loc}(\Omega)$.
  In any case, we can conclude that  $\bar u'\in \cC^{0,\alpha}_{\loc}(\Omega)$, so the proof is complete.
 
   \qed


%
%
%
%

\bigskip

\noindent {\bf Funding:} 
The authors have been supported by the INdAM-GNAMPA group;
R.M. acknowledges also the MIUR Excellence Department Project awarded to the Department of Mathematics, University of Rome ``Tor Vergata'', CUP E83C18000100006.


\bigskip

\noindent {\bf  Compliance with Ethical Standards:}
We wish to confirm that there are no known conflicts of interest
associated with this publication and there has been no significant
financial support for this work that could have influenced its
outcome.



\end{document}